\documentclass[reqno,11pt]{amsart}
\usepackage{amsmath, latexsym, amsfonts, amssymb, amsthm, amscd}
\usepackage{epsfig}
\usepackage{graphics,epsf,psfrag}

\numberwithin{equation}{section}

\newtheorem{theorem}{Theorem}[section]
\newtheorem{lemma}[theorem]{Lemma}
\newtheorem{proposition}[theorem]{Proposition}

\newtheorem{rem}[theorem]{Remark}

\newcommand{\bbE}{{\ensuremath{\mathbb E}} }

\newcommand{\bbN}{{\ensuremath{\mathbb N}} }

\newcommand{\bbP}{{\ensuremath{\mathbb P}} }

\newcommand{\bbS}{{\ensuremath{\mathbb S}} }

\newcommand{\cA}{{\ensuremath{\mathcal A}} }

\newcommand{\cC}{{\ensuremath{\mathcal C}} }
\newcommand{\cD}{{\ensuremath{\mathcal D}} }

\newcommand{\cG}{{\ensuremath{\mathcal G}} }
\newcommand{\cH}{{\ensuremath{\mathcal H}} }

\newcommand{\cL}{{\ensuremath{\mathcal L}} }

\newcommand{\cN}{{\ensuremath{\mathcal N}} }

\newcommand{\cP}{{\ensuremath{\mathcal P}} }

\newcommand{\cT}{{\ensuremath{\mathcal T}} }

\newcommand{\ga}{\alpha}
\newcommand{\gb}{\beta}
\newcommand{\gga}{\gamma}            
\newcommand{\gG}{\Gamma}
\newcommand{\gd}{\delta}

\newcommand{\gs}{\sigma}

\newcommand{\go}{\omega}

\renewcommand{\tilde}{\widetilde}          
\DeclareMathSymbol{\leqslant}{\mathalpha}{AMSa}{"36} 
\DeclareMathSymbol{\geqslant}{\mathalpha}{AMSa}{"3E} 
\DeclareMathSymbol{\eset}{\mathalpha}{AMSb}{"3F}     
\newcommand{\dd}{\text{\rm d}}             

\newcommand{\sumtwo}[2]{\sum_{\substack{#1 \\ #2}}} 
\newcommand{\limtwo}[2]{\lim_{\substack{#1 \\ #2}}}     

\newcommand{\R}{\mathbb{R}}
\newcommand{\Z}{\mathbb{Z}}
\newcommand{\N}{\mathbb{N}}

\DeclareMathOperator{\sign}{sign}

\def\bP{\ensuremath{\bs{\mathrm{P}}}} 
\def\bE{\ensuremath{\bs{\mathrm{E}}}}

\newcommand{\ind}{\bs{1}}

\def\bs{\boldsymbol}

\def\rc{\mathrm{c}}
\def\rf{\mathrm{f}}

\def\tL{\tilde{L}}

\def\tf{\textsc{f}}

\def\tL{\tilde L}
\def\eig{\mathtt{Z}}
\def\gab{{\ga,\gb}}
\def\bQ{\bs{Q}}
\def\bq{\bs{q}}
\def\bp{\bs{\mathrm{p}}}

\newcommand{\obbS}{\overline{\ensuremath{\mathbb S}} }
\newcommand{\obbN}{\overline{\ensuremath{\mathbb N}} }

\title[Infinite volume limits of periodic polymer chains]{Infinite volume limits\\
of polymer chains with periodic charges}

\author{Francesco Caravenna}

\address{
Institut f\"ur Mathematik, Universit\"at Z\"urich,
Winterthurerstrasse 190, CH--8057
Z\"urich}

\email{francesco.caravenna\@@math.unizh.ch}

\author{Giambattista Giacomin}

\address{Laboratoire de Probabilit{\'e}s de P 6\ \& 7 (CNRS U.M.R. 7599) and  Universit{\'e} Paris 7 -- Denis Diderot, U.F.R. Mathematiques, Case 7012, 2 place Jussieu, 75251 Paris cedex 05, France
}

\email{giacomin\@@math.jussieu.fr}

\author{Lorenzo Zambotti}
\address{Dipartimento di Matematica, Politecnico di Milano,
Piazza Leonardo da Vinci 32, 20133 Milano, Italy}
\email{lorenzo.zambotti\@@polimi.it}

\date{\today}

\begin{document}

\begin{abstract}
The aim of this paper is twofold:
\\
-- To give an elementary and self-contained proof of an explicit
formula for the free energy for a general class of polymer chains
interacting with an environment through periodic potentials. This
generalizes a result in \cite{cf:BG} in which the formula is derived
by using the Donsker-Varadhan Large Deviations theory for Markov chains.
We exploit instead tools from renewal theory.
\\
-- To identify the infinite volume limits of the system. 
In particular, in the
different regimes we encounter transient, null recurrent and positive recurrent
processes (which correspond to {\sl delocalized},
{\sl critical} and {\sl localized} behaviors of 
the trajectories). This is done by exploiting the sharp estimates on
the partition function of the system obtained by the renewal theory approach.

The precise characterization of the infinite volume limits of the system
exposes a non-uniqueness problem. We will
however explain in detail how this (at first) surprising phenomenon
is instead due to the presence of a first-order phase transition.
\end{abstract}

\maketitle


\section{Introduction and main results}
\label{sec:intro}

\begin{figure}[h]
\begin{center}
\leavevmode
\epsfxsize =14.0 cm
\psfragscanon
\psfrag{A}[c][c]{A}
\psfrag{B}[c][c]{B}
\psfrag{Af}[c][c]{\footnotesize A--favorable}
\psfrag{Bf}[c][c]{\footnotesize B--favorable}
\psfrag{interface}[c][c]{interface layer}
\epsfbox{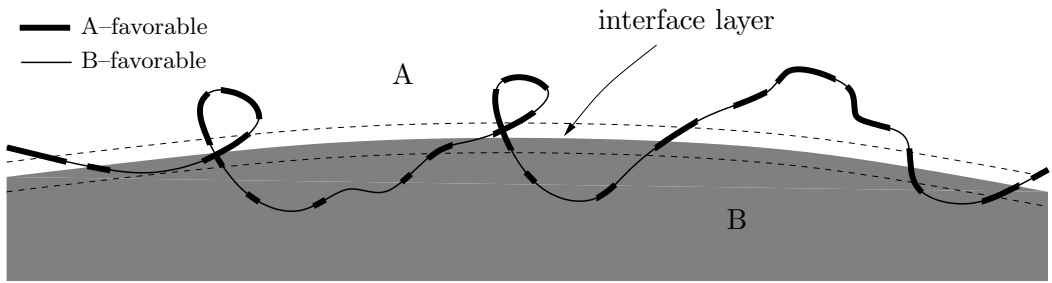}
\end{center}
\caption{\label{fig:cop_ads} \footnotesize
A polymer is made up of two types of monomers, type A that lies preferably
in the solvent A and type B that prefers solvent B (a thicker line
denotes the stretches of 
A--favorable monomers). In order to satisfy
the preferences of all the monomers, the polymer has to keep
close to the interface between the two solvents.
However entropy plays a role too and observing perfect monomer--solvent
  matchings is highly improbable as long as the temperature of the
  system is non zero. A realistic model may also include a different
  type of interaction, attractive or repulsive, at the interface layer.
}
\end{figure}

The real systems that we want to model
are schematized in Fig.~\ref{fig:cop_ads}.
A linear {\sl polymer}, that is a chain made up of
{\sl almost} repetitive units (the {\sl monomers}),
fluctuates in a medium constituted by  two solvents,
A and B, separated by an interface. We say {\sl almost}
repetitive because the monomers differ for one property,
that we call {\sl charge}, that determines the affinity
of the monomer for one or the other solvent (in the figure the charge 
is considered simply as positive, i.e. A--favorable, or negative,
i.e. B--favorable, but in general it may have  an intensity which 
also varies from 
monomer to monomer).

Let us consider the following two possible scenarios:
\smallskip

\begin{itemize}
\item[1--] Imagine that there are as many monomers preferring
the solvent A as the ones preferring B and 
that the charges are distributed along the chain in such a way
that, roughly, the charges alternate. Then the only
configurations with all monomers in their preferred solvent
are configurations that stick closely to the interface. This is
true
even if the matching of charges and solvents is only approximate.
If this is what happens, we say that the polymer is localized
at the interface.
\item [2--] The limit of the argument above is that it takes into account
only of energetic effects (the charge dependent
 interaction monomer--solvent). In particular perfect matchings
 are essentially impossible in large systems with non zero temperature,
 but imperfect matchings leave open the possibility of observing
 the localization phenomenon outlined above.
 In reality however, if for example
 the A--favorable monomers outnumber the B--favorable
 ones, say in a ratio two to one,   then it is still true
 that the polymer may end up optimizing the energetic
 gain via (possibly imperfect) matchings, but it may also take
 the different strategy of lying above the interface
 almost all the time, performing in this way only very imperfect
 matchings, two over three, but gaining (presumably)
 very much in fluctuation freedom (the so called {\sl
 entropic gain}).
 \end{itemize}
\smallskip

The situation is therefore rather unclear
 and it appears that a non trivial energy-entropy competition is
 governing the system.

The type of polymer we have introduced is what is called
a {\sl copolymer}, a synonymous with heterogeneous polymer,
and the physical system 
goes under the name of {\sl copolymer near an interface between
selective solvents} \cite{cf:BdH,cf:MGO,cf:G}.
More realistic would be however to consider that the interface
is typically non extremely sharp and there is a layer
in which the two solvents mix (we could also imagine that in this
layer are trapped some impurities) 
and
the monomers may pay a price or receive a reward in crossing,
or even lying on,
the interface layer. This extra interaction is
normally referred to as a {\sl pinning} or {\sl depinning} interaction
(but also as {\sl adsorption/desorption}) \cite{cf:W}.

One can also imagine the extreme case in which there
is no monomer-solvent interaction, but there are only
(de)pinning interactions: this is a very realistic situation too,
even beyond the two solvent picture we have given.
One can in fact imagine that a polymer fluctuates freely
in space (or in a solvent), except when it is in proximity
of a {\sl defect line} with which it interacts, see e.g. \cite{cf:AS} and
references therein.

The literature on the realistic situations that we have just
outlined is vast. A  considerable part of it focuses
on a case which is very relevant for applications: the one
of periodic distribution of charges (we mention of course
also the other extremely relevant case of
disordered charge distributions \cite{cf:AS,cf:BdH,cf:G}). 
With this we mean that
the sequence of charges repeats after a finite
number of monomer units and the polymer is
effectively made up piecing together identical
stretches of monomers. In this paper we focus exactly on
modeling periodic copolymer models with adsorption
(or pinning) interactions.

A further important remark is that a polymer model
should include the so called {\sl excluded volume interaction},
which leads to self-avoiding walks.
We will enforce  the self-avoiding condition
 by making the rather drastic choice of considering directed polymers
(of course the trajectory in
 Fig.~\ref{fig:cop_ads} may be self-avoiding,
 if the space is three dimensional).

\subsection{The model}
We consider a random walk $S:=\{S_n\}_{n\in \N\cup \{0\}}$
that is $S_0=0$ and $S_n:= \sum_{j=1}^n X_j$, where
$\{X_j\}_{j \in \N}$ is an IID sequence with $\bP( X_1=+1)=
\bP( X_1=-1)=p\in (0,1/2)$ and $\bP( X_1=0)=1-2p$ (we have
decided to exclude $p=1/2$ only for notational convenience,
because of the periodicity of the walk).

The walk $S$ is our free model. We suppose that $S$ interacts with an
environment, that we model with four periodic sequences
$\go^{(+1)}$, $\go^{(-1)}$, $ \go^{(0)}$ and  $\widetilde\go^{(0)}$.
We consider two models, respectively {\sl free} and
{\sl constrained}, defined by 
\begin{equation}
\label{eq:genP}
\frac{\dd \bP^\rf_{N, \go}}{\dd \bP} (S)\, =\,
\frac{\exp\left(\cH_N (S)\right)}{\tilde Z_{N, \go}^\rf}
\ \ \ \text{ and } \ \ \
\frac{\dd \bP^\rc_{N, \go}}{\dd \bP} (S)\, =\,
\frac{\exp\left(\cH_N (S)\right)}{\tilde Z_{N, \go}^\rc}
 \ind_{\left\{ S_N=0 \right\}}\,,
\end{equation}
where the Hamiltonian is 
\begin{equation}
 \label{eq:genH}
\cH _N (S)\, :=\,
 \sum_{i =\pm 1} \sum _{n=1}^ {N}
 \go_n^{(i)} \ind_{\left\{\sign \left( S_n\right)=i\right\}}+
 \sum _{n=1}^ {N}
 \go^{(0)}_n \ind_{\left\{S_n=0 \right\}}+
  \sum _{n=1}^ {N}
 \widetilde\go^{(0)}_n \ind_{\left\{\sign\left(S_n\right)=0 \right\}}\,.
 \end{equation}
Some comments are in order:
\begin{enumerate}
\item
$\go^{(\pm1)}$, $ \go^{(0)}$ and  $\widetilde\go^{(0)}$
are periodic sequences of real numbers, describing the interaction of
the monomers with the solvents and the interface. We say that
the sequence $\go=\{ \go_n\}_{n \in \N}$ is periodic
if there exists $T\in \N$ such that
$\go_{n+T}= \go_n$ for every $n$. The smallest such $T$
is the period of $\go$.
From now on
$\go$ rather denotes  the four periodic sequences appearing in
\eqref{eq:genH}, and we will use $T=T(\go)$ to denote
the least  common multiple
of the periods of $\go^{(\pm1)}$, $ \go^{(0)}$ and  $\widetilde\go^{(0)}$.
\item To define $\sign(S_n)$ when $S_n=0$ we adopt the following convention:
if $S_{n-1} \ne 0$ we set $\sign (S_n):= \sign(S_{n-1})$ while if also $S_{n-1}=0$ we set
$\sign(S_n):=0$. This definition has the following simple interpretation: $\sign(S_n)=+1,-1$ or $0$
according to whether the bond $(S_{n-1},S_n)$ lies above, below or on the $x$--axis.
\item $\tilde Z_{N, \go}^a := \bE \big[\exp (\cH_N) \, \big(\ind_{\{a=\rf\} } +
\ind_{\{a=\rc\} } \ind_{\{ S_N=0 \}}\big) \big]$ is the normalization constant, that is
usually called {\sl partition function}.
\item The measure~$\bP^a_{N,\go}$ is invariant under
the joint transformation $S \to -S$ and
$\go^{(+1)} \to \go^{(-1)}$, hence  we may (and will) assume that
\begin{equation}
\label{eq:ass_h}
h_\go \, := \, \frac 1{T(\go)}\sum_{n=1}^{T(\go)}
\left(\go_n^{(+1)}-\go_n^{(-1)}\right) \, \geq \, 0\,.
\end{equation}
\end{enumerate}

\medskip
\begin{rem}[\it Copolymers and pinning models]
\label{rem:motiv}
\rm
The general model \eqref{eq:genP} that we consider will be referred to
as a {\sl copolymer with adsorption} model.
This includes as special cases the
{\sl copolymer} and {\sl pinning} models that were mentioned informally above.
More precisely, the copolymer model corresponds
to the choice $\go^{(0)}= \tilde \go^{(0)} \equiv 0$
(this formulation generalizes the case considered in \cite{cf:BG}).
If instead we set $\go^{(\pm 1)}\equiv 0$ we are left with
the interactions at the interface, or defect line, and we are dealing
with a pinning model.
We stress that much of the literature on periodic
models, e.g. \cite{cf:MGO,cf:NN,cf:NZ} for the copolymer case
and \cite{cf:GG} for the pinning case,
focuses on the case of $T=2$. We mention as exceptions \cite{cf:vRR}
that deals with the free energy of very particular types of periodic sequences and
\cite{cf:SD1,cf:SD2} treating, in a qualitative and non rigorous fashion,
arbitrary $T$ models
 (see \cite{cf:BG} for more details on the literature).
\end{rem}

\subsection{The free energy and the localization/delocalization alternative}
Getting back to the general model (\ref{eq:genP}),
we observe that from a technical viewpoint
it is convenient to set
\begin{equation} \label{eq:newH}
    \cH^\prime_N (S)
    \, := \, \cH _N (S) \;-\; \sum _{n=1}^ {N}\go_n^{(+1)}\,,
\end{equation}
which just corresponds to $\go_n^{(+1)} \to 0$,
$\go_n^{(-1)} \to (\go_n^{(-1)}-\go_n^{(+1)})$,
$\widetilde\go_n^{(0)} \to (\widetilde\go_n^{(0)}-\go_n^{(+1)})$, and to note that
this new energy yields {\sl the same polymer measures}, namely
\begin{equation}
\label{eq:newP}
\frac{\dd \bP^a_{N, \go}}{\dd \bP} (S)\, =\,
\frac{\exp\left(\cH_N' (S)\right)}{Z_{N, \go}^a} \left(\ind_{\left\{a=\rf\right\} } +
\ind_{\left\{a=\rc\right\} } \ind_{\left\{ S_N=0 \right\}}
\right)\,,
\end{equation}
where the new partition function is just
$Z_{N, \go}^a=\tilde Z_{N, \go}^a \exp(- \sum _{n=1}^ {N}\go_n^{(+1)})$.

It is not difficult to see that $\{ \log Z_{nT, \go}^\rc\}_n$ is a super-additive sequence
and from this to establish the existence of the limit
\begin{equation}
\label{eq:fe}
\tf_\go \, :=\, \lim_{N\to \infty} \frac 1N\log Z_{N, \go}^\rc.
\end{equation}
$\tf_\go$ is the {\sl free energy} of the system.
It is also rather straightforward to show that
\eqref{eq:fe} holds also if we replace the superscript $\rc$ with $\rf$,
i.e. the free energy does not depend on the boundary condition.
For a proof of these facts the reader is referred for example to \cite{cf:G},
but we stress that in this paper we will give a proof
of the existence of the limit in \eqref{eq:fe} that does not rely on
super-additivity, see Section~\ref{sec:fe}. 

Leaving aside for the moment the problem of determining
$\tf_\go$, we focus instead on a simple but crucial aspect
of the free energy, namely that
\begin{equation}
\label{eq:fe>=0}
\tf_\go \, \ge \, 0.
\end{equation}
The proof of this fact is absolutely elementary:
\begin{equation}
\label{eq:proof_fe>=0}
\begin{split}
 \frac 1N\log Z_{N, \go}^\rc\, &\ge \,
 \frac1N \bE\left[ \exp \left(  \cH^\prime_N (S)\right);\; S_n>0 \text{ for } n=1, \ldots, N-1
 \right]
 \\
 & = \frac 1N \log \left(\frac 12 \exp (\go_N^{(0)}) K(N)\right)
 \stackrel{N \to \infty} \longrightarrow 0,
 \end{split}
\end{equation}
where we have introduced the notation $K(N):= \bP (S_n\neq 0 $
for $ n=1, \ldots, N-1$ and $S_N=0)$ and we
have used the polynomial decay of $K(\cdot)$.
Later on we will need the precise asymptotic behavior of $K(\cdot)$.
This can be found in \cite[Ch.~XII]{cf:Feller2} and we anticipate it here:
\begin{equation}
\label{eq:Feller2}
K(N) \, \stackrel{N \to \infty} \sim \, \frac{c_K}{N^{3/2}},
\end{equation}
for a positive constant $c_K$ that depends on $p$ (by $a_N \sim b_N$
we mean of course $a_N/b_N \to 1$ as $N\to\infty$).

\smallskip

Inspired by \eqref{eq:fe>=0} and by its proof,
it is customary to say that the system is
\begin{itemize}
\item[-] localized if $\tf_\go >0$;
\item[-] \rule{0pt}{1.1em}delocalized if $\tf_\go =0$.
\end{itemize}
As unsatisfactory as this definition may look at first,
we will see in the next paragraph that the above {\sl dichotomy} captures
some of the essential features of the system.
For the moment we would like to stress that
the free energy $\tf_\go$ admits an {\sl explicit formula} in terms
of the charges~$\go$, see Theorem~\ref{th:fe} below,
that has been first derived in \cite{cf:BG}, by means of large deviations
techniques, and then re--obtained in \cite{cf:CGZ1},
using a more direct approach based on renewal theory. One of the purposes of
this work is to present (in Section~\ref{sec:fe}) a direct self-contained proof of this formula,
using renewal theory ideas in analogy to~\cite{cf:CGZ1}.
As we shall see next, the renewal theory approach
allows to go much further.

\subsection{From free energy to path behavior}
\label{sec:fetopath}
A very natural question is whether the localization (resp. delocalization)
defined in terms of the free energy does correspond to a real localized
(resp. delocalized) behavior of the trajectories of the polymer measure $\bP_{N,\go}^a$.
A positive answer to this question had been already given before,
but only in terms of weak (de)localization results and leaving out essentially
in all instances the critical behavior (see \cite{cf:BG} and references therein).
We have instead given strong path results in terms of scaling
limits in \cite{cf:CGZ1}, by exploiting renewal theory ideas.
Here we pursue the line and obtain the precise
characterization of the infinite volume limit of the system.

The key technical point is that we can
 go well beyond the Laplace asymptotic
behavior captured by the free energy.
In fact in \cite{cf:CGZ1} we have shown that there exists
a basic parameter $\gd_\go$, which is
an explicit function of the charges~$\go$, that determines the
{\sl precise asymptotic behavior} of the partition function
(we define $\gd_\go$ in \eqref{eq:gdgb}, but the precise expression
of $\gd_\go$ is not essential now). Let us denote by $\bbS$
the Abelian group
$\Z / (T\Z)$, that is $\{0,\ldots,T-1\}$ with sum modulo~$T$,
and we write equivalently $[n]=\ga$ or
$n \in \ga$ to denote that $n$ is in the equivalence class of~$\ga \in \bbS$.
The result proven in \cite{cf:CGZ1} is:
\medskip

\begin{theorem}[\bf Sharp asymptotic estimates] \label{th:sharpest}
Fix $\eta \in \bbS$ and consider the asymptotic behavior of $Z_{N,\go}^\rc$ as
$N\to\infty$ along $[N]=\eta$. Then:
\begin{enumerate}
\item If  $\gd_\go < 1$ then
$Z_{N,\go}^\rc \;\sim\; C^<_{\go,\eta} \, / \, N^{3/2}$ \,;
\item\rule{0pt}{1.1em}If $\gd_\go = 1$ then $Z_{N,\go}^\rc \;\sim\;
C^=_{\go,\eta} \, / \, N^{1/2} $;
\item\rule{0pt}{1.1em}If  $\gd_\go > 1$   then $\tf_\go>0$ and
$Z_{N,\go}^\rc \;\sim\; C^>_{\go,\eta} \, \exp \big(\tf_\go N
\big)$,
\end{enumerate}
where the positive quantities $\tf_\go$, $C^>_{\go,\eta}$, $C^<_{\go,\eta}$ and $C^=_{\go,\eta}$ are
given explicitly in Theorem~\ref{th:as_Z}.
\end{theorem}
\medskip

In Theorem~\ref{th:as_Z} one finds also the asymptotic behavior
for the free endpoint case.
We remark   that Theorem~\ref{th:sharpest} implies that
 the {\sl localized regime}
 corresponds to $\gd_\go > 1$.
The complementary delocalized regime $\gd_\go \le 1$ clearly
splits in two sub-regimes
that we call {\sl strictly delocalized regime} ($\gd_\go < 1$) and
{\sl critical regime} ($\gd_\go = 1$).
The reason for such a denomination is clear if one considers
that $\go \mapsto \gd_\go$ is a continuous function on the set
$\{\go:\; T(\go)=T\}$ (that is for fixed period) and hence arbitrarily
small variations in $\go$ may change $\gd_\go=1$ to
$\gd_\go>1$ or $\gd_\go<1$, while of course the localized
and strictly delocalized regimes are {\sl stable}.

Theorem~\ref{th:sharpest} has been applied in \cite{cf:CGZ1} to determine
the scaling limits of our models. More precisely, it has been shown that
for every fixed $\eta \in \bbS$ the linear interpolation of $\{S_{i/N}/\sqrt{N}\}_{i=0, \ldots, N}$
under $\bP_{N,\go}^a$ converges in distribution as $N\to\infty$ along the subsequence $[N]=\eta$.
The properties of the limit process (that in general may depend on the choice of $\eta$)
are radically different in the three regimes $\gd_\go \lesseqqgtr 1$
and this gives a precise picture of localization/delocalization (see \cite[Th.~1.3]{cf:CGZ1}).

It is natural to look at the scaling limits as describing
the {\sl global properties} of the system.
In this paper we focus rather on the {\sl infinite volume limit} of our model,
that is on the weak convergence of the polymer measures
$\bP^a_{N, \go}$ {\sl without rescaling}, as a measure on $\Z^{\N\cup\{0\}}$.
The latter space being equipped with the product topology, weak convergence simply means 
convergence of all finite dimensional marginal distributions and hence
the infinite volume limit contains the information on the {\sl local properties}
of the model.

In the following theorem, that is our main result,
we characterize the possible limits of $\bP^a_{N, \go}$, showing that they exhibit distinctive
features of localization/delcalization according to whether $\gd_\go > 1$ or $\gd_\go < 1$
(the critical case $\gd_\go = 1$ is  borderline, as for the scaling limits).

\begin{theorem}[\bf Infinite volume limit]\label{th:infvol}
For every $\eta \in \bbS$ and for $a=\rf, \rc$ the polymer measure $\bP^a_{N, \go}$
converges weakly as $N\to\infty$ along the subsequences $[N]=\eta$ to a limit measure
$\bP_{\go}^{\eta, a}$, law of an irreducible Markov process on $\Z$ which is:
\begin{itemize}
\item positive recurrent if $\gd_\go>1$ (\textsl{localized regime})\,;
\item \rule{0pt}{1.1em}null recurrent if $\gd_\go=1$ (\textsl{critical regime})\,;
\item \rule{0pt}{1.1em}transient if $\gd_\go<1$ (\textsl{strictly delocalized regime})\,.
\end{itemize}
When $\gd_\go \ge 1$ the limit 
law $\bP_{\go}^{\eta, a} = \bP_\go$  does not depend
on $\eta$ and $a$, hence both the polymer measures
$\bP^\rf_{N, \go}$ and $\bP^\rc_{N, \go}$ converge weakly as $N\to\infty$ to the same limit $\bP_\go$.
\end{theorem}
\noindent
We prove this theorem in Section~\ref{sec:infvol}, exploiting
the precise asymptotic behavior of $Z_{N,\go}^a$ given in Theorem~\ref{th:sharpest}
and in Proposition~\ref{th:as_Z}, and we also provide an explicit construction
of the limit law $\bP^{\eta,a}_{\go}$ in all regimes. We also
point out that the transition kernel of the Markov law $\bP^{\eta,a}_{\go}$ is
only {\sl periodically inhomogeneous},
that is $\bP^{\eta,a}_{\go}(S_{n+1}=y | S_n=x)$ is a $T$--periodic function of~$n$.

Results similar to Theorems~\ref{th:sharpest} and~\ref{th:infvol} have been
obtained for homogeneous pinning systems (see \cite{cf:CGZ2,cf:DGZ,cf:IY})
and for periodic pinning models in the $T=2$ case \cite{cf:MGO}
(we stress however that a $T=2$ periodic pinning model based on
simple random walk
becomes, by considering the marginal on odd or even sites,
a homogeneous model based on a random walk with jumps in
$\{-1,0,1\}$: this {\sl decimation procedure} is
 less straightforward  for $T>2$ and it leads to rather involved models).

In spite of recent advances, see \cite{cf:GT,cf:GTloc} and references therein,
obtaining results like  Theorem~\ref{th:sharpest} and Theorem~\ref{th:infvol}
for disordered models appears to be a real challenge (the problem is more apparent
for the delocalized regime, but also the localized regime of
disordered systems is still only partly understood).

\subsection{Non-uniqueness and first order transition}
\label{sec:P}

It should have possibly struck the reader the dependence of the infinite volume
limit on the {\sl boundary conditions} $[N]$ and on $a=\rf, \rc$
in the strictly delocalized regime $\gd_\go < 1$ (recall that our system is one dimensional!).
This {\sl trouble} was already present in \cite[Th.~1.3]{cf:CGZ1}, i.e. for the scaling limits,
where however also the critical regime is affected.

Here we are going to clarify this point. First of all
we point out that for a large number of cases, that we characterize explicitly
in \S\ref{sec:tl_del}, all limit laws $\bP^{\eta,a}_{\go}$
appearing in Theorem~\ref{th:infvol} in fact coincide also in the strictly delocalized
regime, hence
there is only one infinite volume measure which is the limit of both $\bP^{\rf}_{N,\go}$
and $\bP^{\rc}_{N,\go}$ as $N\to\infty$. This is true in particular
for copolymer and pinning models (defined in Remark~\ref{rem:motiv}).

However there do exist cases when the laws $\bP^{\eta,a}_{\go}$
have a real dependence on the boundary conditions $a=\rf,\rc$ and $[N]=\eta$
(we anticipate that this happens only for $h_\go=0$).
In Section~\ref{sec:Gibbs} we study in detail this phenomenon,
showing that all possible limit laws $\bP^{\eta,a}_{\go}$
are in fact superpositions of {\sl two extremal Gibbs measures}
$\bQ^+_\go$ and $\bQ^-_\go$,
that we define explicitly and which differ sharply for the asymptotic
behavior as $N\to\infty$: with probability~$1$, $S_N \to + \infty$ under $\bQ^+_\go$
and $S_N \to - \infty$ under $\bQ^-_\go$ (we recall that for $\gd_\go < 1$ the infinite volume
process is transient). We insist however on the fact that, in general,
$\bQ^\pm_\go$ differ also for the statistics of the finitely many returns close to the origin and they are not related by a simple symmetry.

\begin{figure}[h]
\begin{center}
\leavevmode
\epsfysize =7 cm
\psfragscanon
\psfrag{D}[c][c]{\Large $\cD$}
\psfrag{L}[c][c]{ \Large $\cL$}
\psfrag{r1}[c][c]{$\varrho=+1$}
\psfrag{r2}[c][c]{ $\varrho=-1$}
\psfrag{0}[c][c]{$0$}
\psfrag{h}[c][c]{ $h$}
\psfrag{bc}[c][c]{$\gb _c$}
\psfrag{b}[c][c]{$\gb$}
\epsfbox{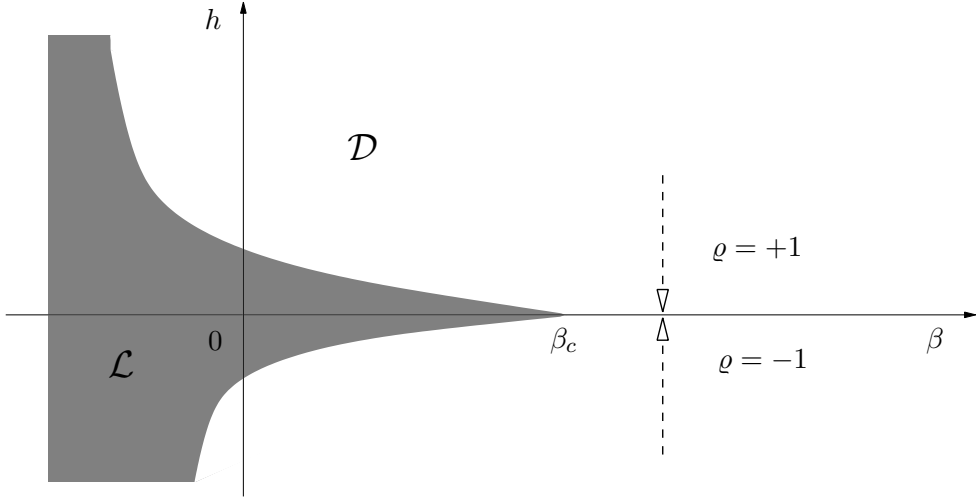}
\end{center}
\caption{\label{fig:P}  \footnotesize
A sketch of the phase diagram for  the model \eqref{eq:Pfig}.
Approaching $h=0$ in the sense of the dashed arrowed lines, one
observes the two sharply different behaviors of paths completely delocalized
above ($\varrho=+1$) or below ($\varrho=-1$) the interface. Taking
the (weak) limits as $h\searrow 0$ (respectively $h\nearrow 0$) of the
infinite volume measures one obtains precisely the measure $\bQ_\go^+$
(respectively $\bQ_\go^-$). The infinite volume limit
for $h=0$ instead exists only along subsequences since there are in general
$T$ different limit points (for the constrained endpoint case and $T$
for the free endpoint case) that are different superpositions of $\bQ_\go^+$
and $\bQ_\go^-$.
}
\end{figure}

We stress that this multiplicity of infinite volume measures should not be regarded
as a pathology, but it is rather the sign of the presence of a
{\sl first order phase transition} in the system. In order to be more precise
let us consider for instance the case of
\begin{equation}
\label{eq:Pfig}
 \frac{\dd \bP_{N, \go}}{\dd \bP} (S) \, = \, \frac{1}{\tilde Z_{N,\go}} \, \exp \left( \sum _{n=1}^ {N}\left(
 \go_n +h \right)\sign \left( S_n\right) -\gb \sum _{n=1}^ {N}
  \ind_{\{ S_n=0\}}\right),
\end{equation}
with $h$ and $\gb$ two real parameters and $\go$ a fixed centered
($\sum_{n=1}^T\go_n=0$) periodic configuration of charges which is non trivial,
that is $\go_i \ne 0$ for some~$i$. For the sake of this paragraph we
define the free energy directly by $f(\gb, h) := \lim_{N\to\infty} N^{-1} \log \tilde Z_{N,\go}$,
that is we do not make the transformation \eqref{eq:newH}. Then with arguments analogous
to \eqref{eq:proof_fe>=0} one gets $f(\gb, h) \ge |h|$ and hence we say that the
system is localized if $f(\gb,h)>|h|$ and delocalized if $f(\gb,h) = |h|$.

The phase diagram of such a model is sketched in Figure~\ref{fig:P}.
In particular
it is easy to show that for $\gb$ sufficiently large and positive
the system is delocalized for any value of~$h$.
On the other hand, for $(\gb, h)=(0,0)$ the system is localized, see \cite[App.~B]{cf:CGZ1} or \cite{cf:BG}.
By the monotonicity of the free energy in $\gb$, one immediately
infers that there exists $\gb_c>0$ such that at $h=0$ localization (resp. delocalization)
prevails for $\gb<\gb_c$ (resp. for $\gb\ge \gb_c$).

The interesting point is that the delocalized regime that appears when $\gb \ge \gb_c$
has sharply different properties according to the sign of~$h$: in fact, since $f(\gb, h) = |h|$,
the quantity $\varrho (\gb,h) := \partial f(\gb,h) / \partial h$ takes the value $+1$
for $h>0$ and $-1$ for $h<0$. Notice that this quantity has the following direct interpretation ($h\neq 0$):
\begin{equation}
    \varrho (\gb,h):= \lim_{N\to \infty}\bE_{N, \go}
\left[ \frac 1N\sum _{n=1}^ {N}  \sign \left( S_n\right) \right]\,.
\end{equation}
Since for $\gb \ge \gb_c$ the free energy is not differentiable at $h=0$,
the system is said to undergo a {\sl first order phase transition}.

It is worth remarking that first order phase transitions are usually associated to multiple infinite
volume limits ({\sl phase coexistence}) like the ones obtained in Theorem~\ref{th:infvol}.
In our case we are able to assert with precision
that $\bQ^\pm_\go$ are {\sl pure phases} (that is extremal Gibbs
states) and 
 which linear combination of
  $\bQ^+_\go$ and $\bQ^-_\go$ one obtains taking the limits along
the subsequences with fixed values of $[N]$.

\subsection{Outline of the paper}
In Section~\ref{sec:fe} we give the formula for $\tf_\go$
and its proof based on renewal theory: along
the proof the fundamental processes characterizing the rest
of the paper will appear naturally. 
Section~\ref{sec:fe} contains only algebraic manipulations
and basic probability facts.
In Section~\ref{sec:she} we recall and discuss a more general version
of Theorem~\ref{th:sharpest}, proven in \cite{cf:CGZ1}.
In Section~\ref{sec:prel} we make a number of manipulations
on the finite volume polymer measures that clarify the role
of the random set of contacts of the polymer with the interface
and the excursions of the polymer in the solvents.
In Section~\ref{sec:infvol} we identify the infinite volume limits of the system,
proving in particular Theorem~\ref{th:infvol}. Finally, in Section~\ref{sec:Gibbs}
we unravel the non-uniqueness phenomenon encountered by taking
infinite volume limits.

\section{A renewal theory path to an explicit expression for the free energy}
\label{sec:fe}

We are going to explain how renewal theory ideas lead to
a representation formula for the partition function $Z_{N,\go}^\rc$ that we
exploit to establish an explicit formula for the free energy $\tf_\go$.
It will be clear that one can go much beyond with such
a formula and we will explain (without a full proof) how to
obtain Theorem~\ref{th:sharpest} from it.

\subsection{The matrix encoding procedure.}
In order to give the formula for the free energy we need to
recall the {\sl matrix encoding procedure}
presented in \cite{cf:BG}.
We recall the definition $\bbS := \Z / T \Z$ and, for $n \in \Z$, we denote
by $[n]\in \bbS$ the equivalence class of $n$,
that is if $m \in [n]$ there exists $j \in \Z$ such that $m=n+ jT$.

The basic structure underlying $S$ is for us the renewal process
$\tau:=\{ \tau_j\}_{j=0,1, \ldots}$ defined
by $\tau_0:=0$ and
\begin{equation}
\label{eq:tau}
\tau_{j+1}\, :=\, \inf \left\{ n> \tau_j:\; S_n=0 \right\},
\end{equation}
and, since $S$ is recurrent, $\tau_j < \infty$ for every $j$, $\bP$--a.s..
The sequence $\tau$, which we will view also as a random
subset of $\N \cup \{0\}$, is a renewal process precisely because
$\{\tau_{j}-\tau_{j-1}\}_{j\in \N}$ is an IID sequence.
It is therefore fully characterized  by the law of $\tau_1$
and we have already set the notation $\bP (\tau_1=n)=K(n)$.
Note that, by \eqref{eq:Feller2}, $S$ is only {\sl null}
 recurrent, since
$\bE[\tau]=+\infty$. In renewal terms, $\tau$ is persistent (but
we will prefer to refer to it as recurrent)
and in fact {\sl null} persistent.

Next we can define a $\bbS \times \bbS$ matrix $\Sigma_{\ga, \gb}$ by the relation
\begin{equation} \label{eq:def_Sigma}
    \sum_{n=n_1+1}^{n_2} (\go_n ^{(-1)}-\go_n ^{(+1)}) \;=\; -(n_2-n_1) \, h_\go
\;+\; \Sigma_{[n_1],[n_2]}\,,
\end{equation}
where $h_\go$ has been defined in \eqref{eq:ass_h} (we stress that the matrix
$\Sigma_{\ga, \gb}$ is well-defined because the charges $\go$ are $T$--periodic).
In this way we have decomposed the above sum into a drift term and a fluctuating term,
where the latter has the key property of depending on $n_1$ and on $n_2$ only
through the respective equivalence
 classes $[n_1]$ and $[n_2]$. Now
 for $\ga, \gb\in \bbS$ and $\ell \in \N $ we define
\begin{equation}
\label{eq:def_Phi}
\Phi^\go_{\ga, \gb}(\ell)\, :=\,
\begin{cases}
    \go ^{(0)}_{\gb} \;+\; \Big(\tilde  \go ^{(0)}_{\gb} - \go^{(+1)}_\gb \Big) & \text{if }
    \ell=1,\ \ell \in \gb-\ga \\
     \go ^{(0)}_{\gb} +
    \rule{0pt}{22pt} \displaystyle  \log  \bigg( \frac 12 \Big(1+
\exp\big( -\ell \, h_\go  + \Sigma_{\ga,\gb}\big)\Big) \bigg) &\text{if } \ell > 1,
\ \ell \in \gb-\ga\\
    \rule{0pt}{16pt} 0 & \text{otherwise}
\end{cases}
\, ;
\end{equation}
and  for $n\in\N$ we introduce the $\bbS \times \bbS$ matrix $M^\go(n)$
defined by
\begin{equation}
\label{eq:matrix}
M^\go_{\ga,\gb}(n) \;:=\; e^{\Phi^\go_{\ga,\gb}(n)}
\, K(n) \, \ind_{(n \in \gb -\ga)}\,.
\end{equation}
Summing over $n\in\N$ the entries of $M^\go$ we obtain a $\bbS \times \bbS$ matrix $B$:
\begin{equation}\label{eq:matrixB}
B_{\ga,\gb} \;:=\; \sum_{n\in\N} M^\go_{\ga,\gb}(n).
\end{equation}
We finally introduce for $b\ge 0$ the $\bbS\times\bbS$ matrix
 $A^\go(b)$:
\begin{equation}
\label{eq:def_A}
    A^\go_{\ga,\gb}(b) := \sum_{n\in\N} M^\go_{\ga,\gb}(n)\exp(-bn)\,.
\end{equation}
Notice that $A^\go(0)=B$. It is important to note that
$A^\go(b)$ is a matrix with positive entries and therefore,
by the classical Perron-Frobenius Theorem \cite{cf:Asm},
its spectral radius $\eig _\go(b)$ is also a positive eigenvalue,
with the property that the corresponding left and right eigenvectors
may be chosen to have strictly positive components. Moreover $\eig_\go(b)$
has also the property of being simple, that is its eigenspace
has dimension one, and it is larger than the  absolute value of
any other (possibly complex) eigenvalue of $A^\go(b)$.

We know also that
$\eig _\go (b)$ is a smooth function of $b$, since $A_\gab^\go(\cdot)$ is smooth
for every $\ga$ and $\gb$, and that $\eig _\go (\cdot)$ is also
strictly decreasing, since the entries $A_\gab^\go(\cdot)$ are.
The inverse function of $\eig_\go(\cdot)$, which is defined on the domain
$(0, \gd_\go]$, will be denoted by $\eig^{-1}_\go(\cdot)$.

We now introduce the basic positive parameter $\gd_\go$, which is defined by
\begin{equation}
\label{eq:gdgb}
\gd_\go \, :=\,  \eig _\go (0)\,.
\end{equation}

\subsection{A matrix representation  and
the   formula for the free energy}
We are now ready to give the explicit formula for the free energy $\tf_\go$:

\medskip
\begin{theorem}
\label{th:fe}
The limit in \eqref{eq:fe} exists and is given by
\begin{equation} \label{eq:fe_formula}
\tf_\go \, =\, \begin{cases}
\eig_\go^{-1} (1) & \mathrm{if} \ \gd_\go > 1 \\
0 & \mathrm{if} \ \gd_\go \le 1
\end{cases}\,.
\end{equation}
\end{theorem}

\medskip

As a preliminary step for the proof of Theorem~\ref{th:fe}
we will make a manipulation on the formula for
$Z_{N,\go}^\rc$ leading to a particularly
useful  matrix expression.
This is in reality very simple, just set $\iota_N:= \sup\{ j:\; \tau_j \le N\}$ and
notice that $\{N \in \tau\} = \{\tau_{\iota_N}=N\}$ is just the event that $\tau_k = N$
for some~$k$. Then, by conditioning on the return times $\tau$ and integrating on the
up--down symmetry of the excursions of~$S$, we can write
\begin{equation}
\label{eq:conv1}
\begin{split}
    Z_{N,\go}^{\rc} & \, = \, \bE \left[
\prod_{j=1}^{\iota_N} \exp\big( \Phi^\go_{[\tau_{j-1}],
[\tau_j]}(\tau_j-\tau_{j-1})\big);
\, N\in\tau
\right] \\
    &\, = \,
    \sum_{k=1}^N \sumtwo{t_0,\ldots,t_k \in \N \cup \{0\}}
{0=:t_0 < t_1 < \ldots < t_k:=N} \, \prod_{j=1}^k
M^\go_{[t_{j-1}], [t_j]} (t_j - t_{j-1}) \,.
\end{split}
\end{equation}
But we can go further with the following algebraic
manipulation: let us denote by $\xi = \xi(b) \in (0, \infty)^\bbS$ the right
eigenvector of $A^\go (b)$ with eigenvalue $\eig_\go(b)$
(the precise normalization is inessential).
Then we introduce the {\sl probability kernel}
\begin{equation}
\label{eq:Gamma}
\gG_\gab (n)\, :=\, \frac 1{\eig_\go(b)} \,
\exp(-bn) \, M^\go_\gab (n) \, \frac{\xi_\gb}{\xi_\ga}\,.
\end{equation}
Equation \eqref{eq:conv1} may then be rewritten as
\begin{equation}
\label{eq:conv2}
    Z_{N,\go}^{\rc}  \, = \, \exp(bN) \frac{\xi_{[0]}}{\xi_{[N]}}
    \sum_{k=1}^N \sumtwo{t_0,\ldots,t_k \in \N \cup \{0\}}
{0=:t_0 < t_1 < \ldots < t_k:=N}
\eig_\go(b)^k
 \prod_{j=1}^k
\gG_{[t_{j-1}], [t_j]} (t_j - t_{j-1}) \,.
\end{equation}

\smallskip

\begin{rem}
\label{rem:Mren} \rm We have called
$\gG$  {\sl probability kernel} because
\begin{equation}
\sum_{n,\gb}\gG_\gab (n)\, =\, \frac{1}{\eig(b)} \, \frac{(A(b)\xi)_\ga}{\xi_\ga}\,=\,1,
\end{equation}
because by definition $\xi$ is the right Perron Frobenius eigenvector
of $A(b)$ (notice that we have dropped the explicit dependence
on $\go$, something that we will frequently do below).
Therefore it is possible to interpret $\Gamma$ as the transition matrix of a Markov
chain on $\bbS\times\N$ that we denote  $\{(J_k, T_k)\}_{k=0,1, \ldots}$:
\begin{equation}
\label{eq:gGtrans}
\bbP_b \left( (J_{k+1}, T_{k+1})=(\gb, n)\big\vert (J_{k}, T_{k})=(\ga, m)
\right)\, =\, \gG_\gab (n)\,.
\end{equation}
Note  that since this transition probability does not depend
on $m$, this chain may be built by first sampling
the $\{J_k\}_k$ process, that is a finite state space ($\bbS$)
Markov chain with transition matrix $\sum_n \gG_\gab (n)$,
and then sampling $\{T_k\}_k$ as independent random
variables with distributions $\gG_{J_{k-1},J_k} (\cdot)/
\sum_n   \gG_{J_{k-1},J_k} (n)$.
\end{rem}

\medskip

Thanks to Remark~\ref{rem:Mren} we interpret \eqref{eq:conv2} in probabilistic terms:
for $J_0:=[0]$ (and $T_0:=0$ for definiteness, but the value of $T_0$ is irrelevant)
we define the {\sl Markov renewal process} $\hat\tau$ as
the partial sum process of the sequence $\{T_k\}_k$, that is
\begin{equation}
\label{eq:Mrenewal}
\hat \tau_j \,:=\,  T_1+\ldots +T_j, \ \ j\in\N, \qquad \quad
\hat \tau_0 := 0 \,.
\end{equation}
This is a particular case of the general
class of Markov renewal processes treated for example in \cite{cf:Asm}.
In terms on this new process, \eqref{eq:conv2} takes a nice
probabilistic expression:

\medskip
\begin{lemma}
\label{th:matrix}
For every $b$ we have
\begin{equation}
\label{eq:conv3}
 Z_{N,\go}^{\rc}  \, = \, \exp(bN) \, \frac{\xi_{[0]}}{\xi_{[N]}} \,
 \bbE_b\left[
 \eig_\go(b)^{\hat\iota_N}; \; N \in \hat\tau
 \right],
\end{equation}
where $\hat\iota_N:= \inf\{j: \; \hat \tau_j \le N\}=
\max (\hat \tau \cap \{0,1, \ldots, N\})$
and we have exploited the fact that $\hat \tau$
may be looked upon as a (random) subset of $\N \cup\{0\}$.
\end{lemma}
\medskip

\noindent
The proof of this lemma follows immediately from \eqref{eq:conv2}, \eqref{eq:gGtrans}
and \eqref{eq:Mrenewal}. Next we pass to the proof of Theorem~\ref{th:fe},
treating separately the three regimes $\gd_\go\gtreqqless 1$.

\medskip

\noindent
{\it Proof of Theorem~\ref{th:fe}, case $\gd_\go>1$.}
Since $\gd_\go = \eig_\go(0) > 1$, the image of $\eig_\go(\cdot)$ contains $1$
and we can set $b:= \eig_\go ^{-1}(1)$, so that $\eig_\go(b)=1$ and \eqref{eq:conv3} becomes
\begin{equation}
\label{eq:conv4}
 Z_{N,\go}^{\rc}  \, = \, \exp(bN) \, \frac{\xi_{[0]}}{\xi_{[N]}} \,
 \bbP_b\left(  N \in \hat\tau \right)\,.
\end{equation}
Since $\xi$ is a vector with positive entries,
\eqref{eq:conv4} implies immediately that the superior limit of
$\{(1/N)
\log Z_{N,\go}^{\rc} \}_N$ is bounded above by $b$
and it suffices to show that $\bbP_b\left( N \in \hat\tau \right)$
does not vanish exponentially fast in $N$
to establish that the free energy exists and that it takes
the value $b$. However it is rather intuitive that
a much better bound holds, namely that there exists $c>0$ such that
\begin{equation}
\label{eq:Ren_th}
\inf_N \bbP_b\left(
  N \in \hat\tau
 \right) \, \ge \, c\,,
\end{equation}
because since $b>0$ the process $\hat \tau$ is positive recurrent,
that is $\sup_j \bbE_b[ T_j]< \infty$.
This is in fact a consequence of the Markov Renewal Theorem \cite[Th.~VII.4.3]{cf:Asm},
which gives the precise asymptotic behavior of $ \bbP_b\left( N \in \hat\tau \right)$
as $N\to\infty$. More directly, it suffices to remark that the processes
$\hat\tau^{\gb}:=\{\hat\tau_j: \; J_j=\gb \}$
is a classical 
(i.e. no Markov dependence) positive recurrent renewal process and that
$\bbP_b\left( N \in \hat\tau \right)= \bbP_b\left( N \in \hat\tau ^\gb \right)$ for $\gb = [N]$.
Therefore the classical Renewal Theorem yields $\bbP_b\left( N \in \hat\tau ^\gb \right) \to
T/\bbE_b[\hat\tau^{\gb}_{2}-\hat\tau^{\gb}_1]>0$ as $N\to\infty$ along the subsequence $[N]=\gb$,
and \eqref{eq:Ren_th} is proven, because there are only finitely many options for $\gb$.
\qed

\medskip

\noindent
{\it Proof of Theorem~\ref{th:fe}, case $\gd_\go=1$.}
Since $\gd_\go = \eig_\go(0) = 1$, also in this case $1$ is in the image of $\eig_\go(\cdot)$
and we set $b= \eig_\go^{-1}(1)=0$. In particular, $\limsup_N(1/N)
\log Z_{N,\go}^{\rc} \le b=0$ like before, but
 we cannot proceed like above for a lower bound,
since, under $\bbP_0$, $\hat \tau$ is null recurrent (that is
$\bbE_0[ \hat\tau_j-\hat \tau_{j-1}]=\infty$).
However, by \eqref{eq:proof_fe>=0}, we already know that
$\liminf_N(1/N)
\log Z_{N,\go}^{\rc} \ge 0$ and we are done.\qed

\medskip

\noindent
{\it Proof of Theorem~\ref{th:fe}, case $\gd_\go<1$.}
This is quick too: since $\gd_\go =\eig_\go(0)<1$, by choosing $b=0$ in
\eqref{eq:conv3} we clearly see  that
$Z_{N, \go}^\rc=O(1)$, so
$\limsup_N(1/N)
\log Z_{N,\go}^{\rc} \le 0$, and
\eqref{eq:proof_fe>=0} provides the lower bound.
Note that in this case to the Markov renewal process
is superimposed a killing rate $\eig_\go(0)$ and it
is this {\sl transient} or {\sl terminating} process
that we should consider as the Markov
renewal process naturally associated to the regime
in which $\eig_\go(0)<1$ (this point will emerge clearly in Section~\ref{sec:infvol}).
\hfill\qed
\medskip

\begin{rem}
\label{rem:split3}\rm
The proof we just completed implicitly contains the
most fundamental ideas of this work, but also
of \cite{cf:CGZ1}. Theorem~\ref{th:sharpest} should appear now
as the natural sharpening of this proof: it is clear that it requires
sharp estimates on suitable {\sl mass renewal functions}, that is
$\bbP_b( N \in \hat \tau)$ and that the three regimes,
corresponding to positive recurrent, null recurrent and
transient Markov renewals, require different techniques
and they in fact yield very different results.
We also stress that the Markov renewal processes arising in the three regimes
are not mere technical tools: they are in fact the limiting
processes given in Section~\ref{th:infvol}.
\end{rem}
\smallskip

\section{Sharp asymptotic estimates}
\label{sec:she}

The aim of this section is to report a more detailed version of Theorem~\ref{th:sharpest},
collecting the results obtained in Section~3 of \cite{cf:CGZ1}, see Theorem~\ref{th:as_Z} below.

We recall from the last section the notation $\xi=\xi(b)$
for the right Perron--Frobenius eigenvector of the matrix $A(b)$,
defined in \eqref{eq:def_A}. More explicitly:
\begin{equation}
    \sum_{\gamma} A_{\ga, \gamma}(b) \, \xi_\gamma \,=\, \eig(b) \, \xi_\ga,
\qquad \forall \ \ga\in\bbS,
\end{equation}
where we recall that $\eig(b)$ is the Perron--Frobenius eigenvalue
of $A(b)$. We choose $\xi$ in $(0,\infty)^\bbS$ and we fix
the normalization $\sum_\gamma \xi_\gamma=1$.
As in the proof of Theorem~\ref{th:fe}, we observe that
when $\gd_\go \ge 1$ the image of $\eig(\cdot)$ contains $1$ and
hence we can set $b := \eig^{-1}(1) = \tf_\go$. From now on, we will
always means that when $\gd_\go \ge 1$ the eigenvector $\xi=\xi(b)$
is evaluated for $b=\tf_\go$ (when $\gd_\go < 1$ we do not
need to use the eigenvector).

\begin{theorem}[\bf Sharp asymptotic estimates] \label{th:as_Z}
Let $k\in\bbN$ with $[k]=\ga$. Then as $N\to\infty$ along $[N]=\eta$ we have:
\begin{enumerate}
\item If  $\gd_\go > 1$ then there exist constants $c^{a}_\eta>0$, $a=\rf,\rc$,
such that:
\begin{equation}\label{eq:as_Z_loc1}
    Z_{N-k,\theta_k\go}^a \; \sim \; \left( c^{a}_\eta\, \xi_\ga \right) \,
\exp \left( \tf_\go\, (N-k) \right)\,.
\end{equation}
\item If $\gd_\go = 1$ then there exist constants $\kappa^{a}_\eta>0$, $a=\rf,\rc$,
such that:
\begin{equation}\label{eq:as_Z_cri}
    Z_{N-k,\theta_k\go}^\rc \;\sim\; \left( \kappa^{\rc}_\eta\, \xi_\ga \right)
\: \frac{1}{\sqrt{N}}\,, \qquad
Z_{N-k,\theta_k\go}^\rf \, \to \, \left( \kappa^{\rf}_\eta\, \xi_\ga \right).
\end{equation}
\item If  $\gd_\go < 1$ then there exist constants $\Lambda_{\ga,\eta}^a>0$, $a=\rf,\rc$,
such that:
\begin{equation}\label{eq:as_Z_del}
    Z_{N-k,\theta_k\go}^\rc  \; \sim \; \Lambda_{\ga,\eta}^\rc
\, \frac{1}{N^{3/2}}\,, \qquad
Z_{N-k,\theta_k\go}^\rf  \; \sim \;  \Lambda_{\ga,\eta}^\rf  \, \frac{1}{N^{1/2}}\,.
\end{equation}
\end{enumerate}
\end{theorem}
The precise value of the constants $\{c^{a}_\eta,\kappa^{a}_\eta,\Lambda_{\ga,\eta}^a\}$
is given in \cite[\S3.2, \S3.3, \S3.4]{cf:CGZ1} and that of $\Lambda_{\ga,\eta}^a$ also in
\eqref{eq:lambda} below. Here we notice that for $\delta_\go\geq 1$
the prefactor in the asymptotic behavior of $Z_{N-k,\theta_k\go}^a$
is equal to a constant, depending on $\eta$ and $a$,
multiplied by the eigenvector $\xi_\ga$: this fact will be important in the proof of Proposition
\ref{pr:infvol1} below.
On the other hand, for $\delta_\go< 1$ in general the constant $\Lambda_{\ga,\eta}^a$
does not admit such a factorization and this is the source of the dependence of the infinite
volume limit on the boundary conditions $a=\rf,\rc$ and $[N]=\eta$.
This phenomenon, anticipated in \S\ref{sec:P}, is studied in detail in Section~\ref{sec:Gibbs}.

\section{The polymer measure: contact set and excursions}
\label{sec:prel}

In this section we perform a preliminary analysis of the polymer measure
$\bP^a_{N,\go}$ that will be a basic tool for the proof of Theorem~\ref{th:infvol},
given in the next section.

The starting point is a very useful decomposition of $\bP^a_{N,\go}$.
The intuitive idea is that a path $\{S_n\}_{n\le N}$ can be split
into two main ingredients:
\begin{itemize}
\item the family $\{\tau_k\}_{k=0,1,\ldots}$ of \textsl{returns to zero}
of $S$, already introduced in \eqref{eq:tau};
\item \rule{0pt}{13pt}the family of \textsl{excursions from zero}
$\{S_{i+\tau_{k-1}}:0\le i\le \tau_{k}-\tau_{k-1}\}_{k=1, 2, \ldots}$
\end{itemize}
Moreover, since each excursion can be
either positive or negative, it is also useful to consider separately
the signs of the excursions $\sigma_k := \sign (S_{\tau_{k-1} + 1})$ and the
absolute values $\{e_k(i):=|S_{i+\tau_{k-1}}|:\, i=1,\ldots,\tau_{k}-\tau_{k-1}\}$.
Observe that these are trivial for an
excursion with length~$1$: in fact if $\tau_{k}=\tau_{k-1} + 1$
then $\gs_k =0$ and~$e_k(0)=e_k(1)=0$.
\smallskip

\begin{rem}\rm
A word about definiteness: if $\tau_k = +\infty$ (and hence $\tau_i = +\infty$
for all $i\ge k$), the definition of the variables $\gs_i$ and $e_i(\cdot)$ given
above do not make sense for $i>k$. However the problem is immaterial, since
in this case these variables are irrelevant
for the purpose of reconstructing the path $\{S_n\}_n$, and consequently we
agree to define $\gs_i$ and $e_i(\cdot)$ for $i>k$ in an arbitrary way.
\end{rem}

\medskip

The process $(\tau_k)_k$ can be also viewed as a (random) subset of $\N \cup \{0\}$,
and for this reason we will refer to it as to the {\sl contact set} (of course we have in mind
the polymer interpretation of our model described in the introduction).
The crucial point, already exploited in \cite{cf:CGZ1} to obtain the scaling
limits our our model, is the following description of the law of the contact set
and of the excursions under the polymer measure $\bP_{N,\go}^a$.

\subsection{The contact set}
\label{sec:zero_set}
We recall the definition $\iota_N = \sup\{k:\, \tau_k\le N\}$.
Let us first consider the returns $(\tau_k)_{k}$ under~$\bP_{N,\go}^a$.
The law of this process can be viewed as a probability measure $\bp^a_{N,\go}$
on the class $\cA_N$ of subsets of $\{1,\ldots, N\}$: indeed
for $A\in\cA_N$, writing
\begin{equation}\label{notA}
A  = \{t_1,\ldots,t_{|A|}\}, \qquad
0 \, =: \, t_0<t_1<\cdots<t_{|A|} \, \leq \, N,
\end{equation}
we can set
\begin{equation}\label{defp^a}
\bp^a_{N,\go}(A) \, := \, \bP^a_{N,\go}(\tau_i= t_i, \ i\leq\iota_N).
\end{equation}
The measure $\bp^a_{N,\go}$ describes the set of contacts of the polymer
with the interface. From the inclusion of $\cA_N$ into $\{0,1\}^{\N\cup\{0\}}$,
the family of all subsets of~$\N\cup\{0\}$, $\bp_{N,\go}^a$ can be viewed as a measure
on~$\{0,1\}^{\N\cup\{0\}}$ (this observation will be useful in the following).

Let us describe more explicitly $\bp^a_{N,\go}(A)$, using the (strong)
Markov property of $\bP^a_{N,\go}$. We use throughout the paper the notation
(\ref{notA}). Recalling the definition \eqref{eq:matrix}
of $M_{\ga,\gb}(t)$, we have for $a=\rc,\rf$:
\begin{equation}\label{eq:relpmz}
\begin{split}
\bp^a_{N,\go} \big( \{k_0, \ldots, k_n\} \big) & \;=\;
\bP^a_{N, \go} \left( \tau_1 =k_1,\ldots,\tau_n=k_n  \right)\\
& \; = \; \left[ \prod_{i=1}^n M_{[k_{i-1}],[k_i]}(k_i-k_{i-1}) \right]
\frac{Z^a_{N-k_n,\theta_{k_n}\go}}{Z^a_{N,\go}},
\end{split}
\end{equation}
for all $0=:k_0<k_1< \cdots < k_n\leq N$ and $a=\rc,\rf$.

\subsection{The signs}
From the very definition \eqref{eq:newP} of our model
it is easy to check that, conditionally on $\{\iota_N, (\tau_j)_{j\le\iota_N}\}$,
the signs $(\sigma_k)_{k\le\iota_N}$ are under $\bP^a_{N,\go}$
an independent family. For $k\le\iota_N$, the conditional law of $\gs_k$ is specified by:
\begin{itemize}
\item[-] if~$\tau_{k} = 1+\tau_{k-1} $, then $\gs_k = 0$;
\item[-] \rule{0pt}{12pt}if~$\tau_{k} > 1+\tau_{k-1} $, then $\gs_k$ can take the two values~$\pm 1$ with
\begin{equation}\label{defsigma}
\bP^a_{N,\go}\Big(\gs_k=+1 \, \Big| \ (\tau_j)_{j\leq\iota_N}\Big)
\, = \, \frac{1}
{1 + \exp\left\{ -(\tau_{k} - \tau_{k-1}) \, h_\go + \Sigma_{[\tau_{k-1}],[\tau_{k}]} \right\}}\,.
\end{equation}
\end{itemize}
Observe that when~$\tau_{\iota_N} < N$ (which can happen only for~$a=\rf$)
there is a last (incomplete) excursion in the interval
$\{0, \ldots, N\}$, and the sign of this excursion is also expressed by~\eqref{defsigma} for~$k=\iota_N+1$,
provided we set~$\tau_{\iota_{N} + 1}:= N$.

\subsection{The moduli of the excursions}
Again, from the definition of our model it follows that,
conditionally on $\{\iota_N, \, (\tau_j)_{j\leq\iota_N},\, (\sigma_j)_{j\leq\iota_N+1}\}$,
the excursions $\big( e_k(\cdot)\big)_{k \le \iota_N+1}$ are under $\bP^a_{N,\go}$ an independent family.
For $k\leq\iota_N$, the conditional law of~$e_k(\cdot)$ on the event
$\{\tau_{k-1}=\ell_0, \ \tau_{k}=\ell_1\}$ is specified for $f=(f_i)_{i=0, \ldots, \ell_1 - \ell_0}$
by
\begin{equation}
\label{defexc}
\begin{split}
& \bP^a_{N,\go}\Big( e_k(\cdot)=f \ \Big| \
\iota_N, \, (\tau_j)_{j\leq\iota_N},\, (\sigma_j)_{j\leq\iota_N+1} \Big) \\
& = \ \bP\Big(S_i=f_i: \ i=0,\ldots,\ell_1-\ell_0 \ \Big| \
S_i > 0: \ i=1, \ldots, \ell_1-\ell_0-1, \ S_{\ell_1-\ell_0}=0\Big)\,.
\end{split}
\end{equation}
For $a=\rf$, when~$\tau_{\iota_N} < N$ the conditional law on the event $\{\tau_{\iota_N}=\ell<N\}$
of the last incomplete excursion $e_{\iota_N + 1}(\cdot)$
is specified for $f=(f_i)_{i=0, \ldots, N - \ell}$ by
\begin{equation}
\label{defexcl}
\begin{split}
&\bP^a_{N,\go}\Big( e_{\iota_N + 1}(\cdot)=f \ \Big| \
\iota_N, \, (\tau_j)_{j\leq\iota_N},\, (\sigma_j)_{j\leq\iota_N+1}
\Big)\\
& = \ \bP\Big(S_i=f_i: \ i=0,\ldots,N-\ell \ \Big| \
S_i> 0: \ i = 1, \ldots, N-\ell\Big).
\end{split}
\end{equation}

\subsection{Building the infinite volume measure}
\label{sec:building}

We stress that the above descriptions of the contact set, of the signs
and of the moduli of the excursions fully characterize
the polymer measure $\bP^a_{N,\go}$. A remarkable fact is that,
conditionally on $(\tau_k)_{k\ge 0}$, the joint distribution
of $(\gs_j,e_j)_{j\leq \iota_N}$ \textsl{does not depend on $N$}: in this sense,
the $N$--dependence is contained in the contact set law $\bp^a_{N,\go}$.

For this reason, the next section is devoted to the study of the
asymptotic behavior of the contact set measure $\bp^a_{N,\go}$ as $N\to\infty$.
The main result is that, for every $\eta \in \bbS$, the measure $\bp^a_{N,\go}$ converges weakly
on~$\{0,1\}^{\N\cup\{0\}}$, as~$N\to\infty$ along the subsequence $[N]=\eta$,
toward a limit measure $\bp^{\eta,a}_\go$ (which in general depends on $a$ and $\eta$).

From this result and from the above considerations, one would like to infer that
the full polymer measure $\bP_{N,\go}^a$ converges weakly on~$\Z^{\N\cup\{0\}}$,
as~$N\to\infty$ along $[N]=\eta$, toward a limit measure $\bP^{\eta,a}_\go$ which is constructed
by {\sl pasting the excursion over the limit contact set}. This is indeed true
when the cardinality of the contact set $\{\tau_n\}_n$ is infinite
under the limit contact set law $\bp_\go := \bp^{\eta,a}_\go$, that is when
$\bp_\go (\tau_k < +\infty) = 1$ for all $k\ge 0$
(we will see that this is what happens when $\gd_\go \ge 1$). In this case the infinite volume
polymer measure $\bP_\go := \bP^{\eta,a}_\go$ can be completely reconstructed from $\bp_\go$
(to lighten the notation, for the rest of this section the dependence
of $\bP^{\eta,a}_\go$ and $\bp^{\eta,a}_\go$ on $a$ and $\eta$ will be omitted).

However when $\gd_\go < 1$ it turns out that the cardinality of the contact set is
$\bp_\go$--a.s. {\sl finite}, hence there is a last infinite excursion.
In this case to obtain the weak convergence of the full polymer measure $\bP_{N,\go}^a$
it is also necessary to determine the law of the sign of the last infinite excursion.
But let us describe more in detail how to construct the infinite volume polymer measure~$\bP_\go$.

\subsubsection{The proper case}
We consider first the case when $\bp_\go (\tau_k < +\infty) = 1$ for all
$k\ge 0$.
Then the infinite volume polymer measure $\bP_\go$ is the law on $\Z^{\N\cup\{0\}}$ under which the processes
$(\tau_j)_j$, $(\sigma_j)_j$ and~$(e_j(\cdot))_j$ have the following laws:
\begin{itemize}
\item The process $(\tau_j)_{j}$ is drawn according to $\bp_\go$.
\item \rule{0pt}{1.1em}Conditionally on $(\tau_j)_{j}$,
the variables $(\sigma_j)_{j}$ are independent.
The conditional law of $\gs_k$ depends only on $(\tau_{k-1},\tau_k)$ and it is specified in the
following way:
\begin{itemize}
\item if $\tau_k - \tau_{k-1} = 1$, then $\gs_k = 0$;
\item if $\tau_k - \tau_{k-1} > 1$, then $\gs_k$ takes the two values $\pm 1$
with probabilities given by the r.h.s. of~\eqref{defsigma}.
\end{itemize}
\item \rule{0pt}{1.1em}Conditionally on $(\tau_j,\gs_j)_{j}$, the variables $(e_j(\cdot))_{j}$
are independent. The conditional law of $e_k(\cdot)$
on the event $\{\tau_{k-1}=\ell_0, \tau_{k}=\ell_1\}$ is given by the r.h.s. of (\ref{defexc}).
\end{itemize}
Of course these requirements determine uniquely the law $\bP_\go$.

\subsubsection{The defective case}

Next we analyze the defective case, when the cardinality of the set
$\{\tau_n\}_n$ is $\bp_\go$--a.s. finite, which is what happens when $\gd_\go < 1$.

Let us denote by $\rho := \sup\{k \ge 0: \tau_k < +\infty\}$ the index of the last point
in the contact set, and by assumption we have $\bp_\go(\rho < +\infty)=1$.
In this case to characterize the infinite volume polymer measure $\bP_\go$
it suffice to specify the laws of the processes $(\tau_j)_{j\in\N\cup\{0\}}$,
$(\sigma_j)_{j = 1, \ldots, \rho + 1}$ and $(e_j(\cdot))_{j = 1, \ldots, \rho + 1}$
under $\bP_\go$.

As before, the process $(\tau_j)_{j}$ is drawn according to the law $\bp_\go$.
Conditionally on $(\tau_j)_{j}$ the variables $(\sigma_j)_{j=1, \ldots, \rho+1}$
are independent, and conditionally on $(\tau_j,\gs_j)_{j}$, the variables
$(e_j(\cdot))_{j=1, \ldots, \rho+1}$ are independent: therefore it remains
to specify the conditional laws of $\gs_k$ and of $e_k(\cdot)$, for $k=1, \ldots, \rho+1$.
However it is easy to see that for $k\le \rho$ there is still no change with respect
to the proper case, that is the conditional laws are given by the r.h.s. of \eqref{defsigma}
and \eqref{defexc} respectively. Hence we are left with specifying the conditional laws
of the last sign $\gs_{\rho + 1}$ and of the last modulus $e_{\rho + 1}(\cdot)$.

For the last modulus the answer is rather intuitive:
on the event $\tau_{\rho + 1} = \ell$, the conditional law of $e_{\rho + 1}(\cdot)$
is given for any~$n\in\N$ and for $f=(f_i)_{i=0, \ldots, n}$ by:
\begin{equation}
\label{defexclul}
\begin{split}
&\bP_{\go}\Big(e_{k}(i)=f_i\,:\ i = 0, \ldots, n \ \Big| \ (\tau_j,\gs_j)_{j} \Big) \, = \,
\bP^\uparrow\Big(S_i=f_i\,:\ i = 0, \ldots, n\Big)
\\ & \qquad := \ \lim_{N\to\infty} \bP\Big(S_i=f_i\,: \ i=0,\ldots,n \ \Big| \
S_i> 0: \ i = 1, \ldots, N\Big),
\end{split}
\end{equation}
where the existence of such limit is well known, cf. \cite{cf:BerDon}.

On the other hand, 
the law of the sign of the last excursion $\gs_{\rho+1}$ has to be
determined by a direct 
computation and this will be done in \S\ref{sec:tl_del}. Once
this is done, the construction of the measure $\bP_\go$ in the defective
case is complete.
A remarkable fact is that, for the choice of free or constrained boundary
conditions, the law of $\gs_{\rho+1}$ is in fact
determined by $\bp_\go$. However  this is not true in
general: one can show (we will  not pursue this point 
in detail) that more general boundary conditions may yield different infinite volume
measures, having the same law for the contact set but a different law for
the sign of the last infinite excursion.

\section{Infinite volume limits}
\label{sec:infvol}

This section contains the proof of Theorem~\ref{th:infvol}. We
study the limit as $N\to\infty$ of the polymer
measure $\bP^a_{N, \go}$, using the sharp asymptotic behavior of the partition
function given in Theorem~\ref{th:as_Z}.
We recall that  $\bP^a_{N,\go}$ is a probability
measure on $\Z^{\N\cup\{0\}}$ and that we endow the latter space with the product topology,
hence weak convergence means convergence of all finite dimensional marginal distributions.

Our focus is mainly on the contact set law $\bp^a_{N,\go}$, defined in \eqref{defp^a},
which is a measure on $\{0,1\}^{\N\cup\{0\}}$. We are going to show that,
for $a=\rf,\rc$, for any fixed $\eta \in \bbS$ and for any value of~$\gd_\go$, the
measure $\bp^a_{N,\go}$ converges weakly on $\{0,1\}^{\N\cup\{0\}}$ as~$N\to\infty$ along the subsequence
$[N]=\eta$. When $\gd_\go \ge 1$ the convergence actually holds true without having to impose
the $[N]=\eta$ constraint, while when $\gd_\go < 1$ the limit may really depend on the
value of~$\eta$ and of $a=\rf,\rc$ (in \S\ref{sec:tl_del} we characterize
precisely the instances in which this happens).

Once the convergence of $\bp^a_{N,\go}$ (as $N\to\infty$ along $[N]=\eta$) is proven,
the analogous statement for the polymer measure $\bP_{N,\go}^a$ follows by
the arguments given in \S\ref{sec:building}.
\smallskip

\begin{rem} \rm
In the proof we actually show that, under the limit measure of $\bp^a_{N,\go}$,
the process $\{\tau_k\}_{k\ge 0}$ is a \textsl{Markov renewal process} with {\sl modulating
chain} $\{J_k\}_{k\ge 0}:= \{[\tau_k]\}_{k\ge 0}$.
This means that, setting $T_k := \tau_k - \tau_{k-1}$ for $k \in \N$,
the joint process $\{(J_k,T_k)\}_{k\in\N}$ is a Markov chain on $\bbS \times \N$
such that the transition probability to go from $(J_{k},T_k)$ to $(J_{k+1},T_{k+1})$ does not
depend on $J_k$:
\begin{equation}
\bbP \left( (J_{k+1}, T_{k+1})=(\gb, n)\big\vert (J_{k}, T_{k})=(\ga, m)
\right)\, =\, \gG_\gab (n)\,,
\end{equation}
see also Remark~\ref{rem:Mren} and the lines that follow it.
The transition kernel $\gG_\gab (n)$ is called the {\sl semi--Markov kernel} of the
Markov renewal process $\{\tau_k\}$. We are going to find an explicit expression
for $\gG_\gab (n)$, showing in particular that the laws of the $T_k$ are:
\begin{enumerate}
\item integrable if $\delta_\go>1$ (localized regime);
\item \rule{0pt}{13pt}defective if $\delta_\go<1$ (strictly delocalized regime);
\item \rule{0pt}{13pt}non integrable if $\delta_\go=1$ (critical regime).
\end{enumerate}
A detailed account on Markov renewal processes can be found in \cite{cf:Asm}.
\end{rem}
\medskip

Next we pass to the proof of Theorem~\ref{th:infvol}. For ease of exposition, we
consider first the cases $\gd_\go > 1$ and $\gd_\go = 1$, where there are no problems of
uniqueness, and then the more delicate strictly delocalized regime $\gd_\go < 1$.

\subsection{The regimes $(\delta_\go>1)$ and $(\delta_\go=1)$}

We are going to prove the following:
\begin{proposition}\label{pr:infvol1}
If $\delta_\go\geq 1$ then
the polymer measures $\bP^\rf_{N, \go}$ and $\bP^\rc_{N, \go}$
converge as $N\to\infty$ to the same limit
$\bP_\go$, under which $(\tau_k)_{k\ge 0}$
is a Markov renewal process with semi--Markov kernel
$(\Gamma_{\ga,\gb}(x): \ga,\gb\in\bbS, x\in\N)$, defined by:
\begin{equation} \label{eq:def_Gamma>}
    \Gamma_{\ga,\gb} (x) \, := \,
M_{\ga,\gb} (x)\, e^{-\tf_\go x}\, \frac{\xi_\gb}{\xi_\ga}.
\end{equation}
\end{proposition}
\noindent
We recall that $\tf_\go=0$ if $\delta_\go\leq1$ and $\tf_\go>0$ if $\delta_\go>1$.

\medskip
\noindent
{\it  Proof of Proposition \ref{pr:infvol1}.}
By the asymptotic behavior of $Z^a_{N,\go}$ in (\ref{eq:as_Z_loc1}) and
(\ref{eq:as_Z_cri}) above,
we have for all $\ga,\eta\in\bbS$ and $k\in\ga$:
\begin{equation}
\exists \limtwo{N\to\infty}{N\in\eta} \
\frac{Z^a_{N-k,\theta_k\go}}{Z^a_{N,\go}}
\, = \, e^{-\tf_\go k} \, \frac{\xi_{[k]}}{\xi_{[0]}},
\end{equation}
and since the right hand side does not depend on $\eta$, then the
limit exists as $N\to\infty$.

By (\ref{eq:relpmz}) it follows that for $0=:k_0<k_1< \cdots < k_j$, $a=\rc,\rf$:
\begin{equation}
\begin{split}
\lim_{N\to \infty}
\bp^a_{N, \go} \big( \{k_0, \ldots, k_j \} \big) & \, =\, \left[
\prod_{i=1}^j M_{[k_{i-1}],[k_i]}(k_i-k_{i-1}) \right]
e^{-\tf_\go k_j} \, \frac{\xi_{[k_j]}}{\xi_{[0]}}
\\ & \, = \, \prod_{i=1}^j \Gamma_{[k_{i-1}],[k_i]}(k_i-k_{i-1})\,,
\end{split}
\end{equation}
and this shows that $\bp^a_{N, \go}$ converges weakly on $\{0,1\}^{\N\cup\{0\}}$
as $N\to\infty$ toward the law $\bp_\go$ under which $(\tau_k)_{k\ge 0}$ is a Markov
renewal process with semi--Markov kernel $\Gamma_{\ga,\gb}(x)$.

Notice that $\sum_{x\in\N,\gb\in\bbS} \Gamma_{\ga,\gb}(x)=1$, that is
$\bp_\go(\tau_k < +\infty)=1$ for all $k \ge 0$. Therefore the weak
convergence of the full polymer measure $\bP^a_{N,\go}$ on $\Z^{\N\cup\{0\}}$ follows
from the arguments given in \S\ref{sec:building}, and the proof is completed.\qed

\subsection{The regime $(\delta_\go<1)$}
\label{sec:tl_del}
We introduce the subset of $\go$ defined by
\begin{equation}\label{sit1}
\cP^< \, := \, \left\{ \go \, : \ \gd_\go < 1, \quad
h_\go=0, \quad \exists \ \ga,\gb: \ \Sigma_{\ga,\gb} \ne 0
\right\}\,,
\end{equation}
where $h_\go$ and $\Sigma_{\ga,\gb}$ have been defined
respectively in \eqref{eq:ass_h} and \eqref{eq:def_Sigma}.
We are going to prove that when $\go \not \in \cP^<$ both the free and
the constrained polymer measures $\bP_{N,\go}^a$, $a=\rf, \rc$, converge
weakly as $N\to\infty$, without having to impose the constraint $[N]=\eta$,
while for $\go \in \cP^<$ the limit exists as $N\to\infty$ along $[N]=\eta$
and in general depend on the choice of $a$
and $\eta$. It is worth stressing that for the two motivating models introduced
in Remark~\ref{rem:motiv}, the pinning and the copolymer models, $\go$ {\sl never}
belongs to $\cP^<$. This is clear for the pinning case,
where by definition $\Sigma\equiv 0$. On the other hand, in the copolymer case
it is known that if $h_\go=0 $ and $\exists \ \ga,\gb: \ \Sigma_{\ga,\gb} \ne 0$
then $\delta_\go>1$, cf. \cite[App.~B]{cf:CGZ1}.

\smallskip

It will turn out that in the strictly delocalized regime
there exists a.s. a last return to zero, i.e. the process  $(\tau_k)_{k\ge 0}$ is defective.
In order to express this with the language of Markov renewal processes,
we introduce the sets
$\obbS:=\bbS\cup\{\infty\}$ and $\obbN:=\bbN\cup\{\infty\}$,
extending the equivalence relation to $\obbN$ by $[\infty] =\infty$.

\smallskip

We need some notation: we set
\begin{equation} \label{eq:as_tPhi}
    \tL_{\ga,\gb} \;:=\;
\begin{cases}
\displaystyle  c_K \big( 1 + \exp(\Sigma_{\ga,\gb}) \big) & \text{if\ \ } h_\go = 0\\
\displaystyle\rule{0pt}{22pt} c_K & \text{if\ \ } h_\go > 0
\end{cases}
\;, \qquad L_{\ga, \gb} \, := \, \frac 12 \, \exp(\go ^{(0)}_{\gb}) \, \tL_{\ga,\gb}.
\end{equation}
We notice that for any $\go$:
\begin{equation}\label{eq:as_L_div}
 L_{\ga, \gb} \, := \, \limtwo{x\to\infty}{[x]=\gb-\ga} \, x^{3/2} \, M_{\ga,\gb}(x).
\end{equation}
In \cite[\S3.4]{cf:CGZ1} it is proven that the constants $\Lambda_{\ga,\eta}^a$
appearing in \eqref{eq:as_Z_del} are equal to:
\begin{equation}\label{eq:lambda}
\Lambda_{\ga,\eta}^\rc \, = \, \big[(1-B)^{-1}L\, (1-B)^{-1}\big]_{\ga,\eta}, \qquad
\Lambda_{\ga,\eta}^\rf \, = \, \big[(1-B)^{-1}\tL\big]_{\ga,\eta},
\end{equation}
where $B$ is defined in \eqref{eq:matrixB}.
Finally we set for all $\ga,\eta\in\bbS$:
\begin{equation}
\mu_{\ga,\eta}^\rc \, := \, \big[L\, (1-B)^{-1}\big]_{\ga,\eta},
\qquad \mu_{\ga,\eta}^\rf \, := \, \tL_{\ga,\eta},
\end{equation}
and for all $\eta\in\bbS$ and $a=\rf,\rc$
we introduce the semi-Markov kernel on $\obbS\times\obbN$:
\begin{equation}\label{eq:GGG}
    \Gamma_{\ga,\gb}^{\eta,a}(x) \, := \,
\begin{cases}
    {\displaystyle M_{\alpha,\gb}(k) \, \Lambda_{\gb,\eta}^a/\Lambda_{\ga,\eta}^a}
\quad & \ga\in\bbS, \ x\in\bbN, \ \beta=[x]\in\bbS\\
    \rule{0pt}{15pt}{\displaystyle \mu_{\ga,\eta}^a/\Lambda_{\ga,\eta}^a} \quad
& \ga\in\bbS, \ x=\infty, \ \beta=[\infty] \\
    \rule{0pt}{15pt}1 & \ga=\gb=[\infty], \ x=0 \\
    \rule{0pt}{15pt}0 & {\rm otherwise}.
\end{cases}
\end{equation}
Notice that $\Gamma^{\eta,a}$ is really a semi-Markov kernel, since for $\ga\in\bbS$:
\begin{eqnarray*}
\sum_{\gb\in\obbS} \sum_{x\in\obbN} \Gamma^{\eta,a}_{\ga,\gb}(x)
& = & \frac{\mu_{\ga,\eta}^a}{\Lambda_{\ga,\eta}^a} +
\sum_{\gb\in\bbS} \sum_{x\in\bbN}
\frac{M_{\alpha,\gb}(x) \, \Lambda_{\gb,\eta}^a}{\Lambda_{\ga,\eta}^a}
= \frac{\mu_{\ga,\eta}^a}{\Lambda_{\ga,\eta}^a} + \frac 1{\Lambda_{\ga,\eta}^a}
[B\cdot \Lambda^a]_{\ga,\eta}
\\ \\ & = & \frac{\mu_{\ga,\eta}^a}{\Lambda_{\ga,\eta}^a} +
\frac 1{\Lambda_{\ga,\eta}^a} (\Lambda_{\ga,\eta}^a-{\mu_{\ga,\eta}^a}) \, = \, 1.
\end{eqnarray*}

\smallskip

We are going to prove the following:

\medskip

\begin{proposition}\label{pr:infvol2}
Let $\delta_\go<1$ and $\eta \in \bbS$. Then:
\begin{enumerate}
\item for $a=\rf,\rc$,
$\bP^a_{N, \go}$
converges weakly as $N\to\infty$ along $[N]=\eta$ toward a measure
$\bP_{\go}^{\eta,a}$, under which $(\tau_k)_{k\ge 0}$
is a Markov renewal process with semi-Markov kernel given by
$\Gamma^{\eta,a}_{\ga,\gb}(x)$.
\item
if $\go\notin\cP^<$, then $\bP_{\go}^{\eta,a}=:\bP_\go$ and $\Gamma^{\eta,a}=:\Gamma^<$
depend neither on $\eta$ nor on $a$, and
both $\bP^\rf_{N, \go}$ and $\bP^\rc_{N, \go}$
converge as $N\to\infty$ to $\bP_\go$, under which $(\tau_k)_{k\ge 0}$
is a Markov renewal process with semi-Markov kernel~$\Gamma^<$.
\end{enumerate}
\end{proposition}

\medskip

\begin{rem}\label{sit2}\rm
Part (2) of Proposition~\ref{pr:infvol2} is an easy consequence of part~(1). In fact
from equation \eqref{eq:as_tPhi} it follows immediately that
when $\go\notin\cP^<$ then both matrices $(L_{\ga,\gb})$ and
$(\tL_{\ga,\gb})$ are constant in $\ga$, and
therefore $\Lambda^a$ factorizes into a tensor product, i.e.
\begin{equation}
\Lambda_{\ga,\eta}^a \, = \, \lambda_\ga^a \, \nu_\eta^a, \qquad
\ga,\eta\in\bbS,
\end{equation}
where $(\lambda_\ga^a)_{\ga\in\bbS}$ and $(\nu_\ga^a)_{\ga\in\bbS}$
are easily computed.
But then it is immediate to check that the semi--Markov kernel
$\Gamma^{\eta,a}=:\Gamma^<$ depends
neither on $\eta$ nor on $a$.
\end{rem}

\bigskip
\noindent
{\it Proof of Proposition \ref{pr:infvol2}.}
By the preceding Remark it suffices to prove part~(1).
By \eqref{eq:as_Z_del} we have
we have for all $\ga,\eta\in\bbS$ and $k\in\ga$:
\begin{equation} \label{eq:miau}
\exists \limtwo{N\to\infty}{N\in\eta} \
\frac{Z^a_{N-k,\theta_k\go}}{Z^a_{N,\go}}
\, = \, \frac{\Lambda_{[k],\eta}^a}{\Lambda_{[0],\eta}^a}.
\end{equation}
By (\ref{eq:relpmz}) it follows that for $0=:k_0<k_1< \cdots <
k_j < \infty$, $a=\rc,\rf$:
\begin{equation}
\begin{split}
\limtwo{N\to\infty}{N\in\eta} \
\bp^a_{N, \go} \big( \{k_0, \ldots, k_j \}
\big) & \, =\, \left[
\prod_{i=1}^j M_{[k_{i-1}],[k_i]}(k_i-k_{i-1}) \right]
\frac{\Lambda_{[k_j],\eta}^a}{\Lambda_{[0],\eta}^a}
\\ & \, = \, \prod_{i=1}^j
\Gamma^{\eta,a}_{[k_{i-1}],[k_i]}(k_i-k_{i-1})\,.
\end{split}
\end{equation}
This shows that $\bp^a_{N, \go}$ converges weakly on $\{0,1\}^{\N\cup\{0\}}$
as $N\to\infty$, $[N]=\eta$, toward the law $\bp^{\eta,a}_\go$
under which $(\tau_k)_{k\ge 0}$ is a Markov
renewal process with semi--Markov kernel $\Gamma^{\eta,a}_{\ga,\gb}(x)$.

However this time the semi--Markov kernel is defective, that is
$\sum_{\gb\in\bbS,x\in\N}\Gamma^{\eta,a}_{\ga,\gb}(x) < 1$, hence the
contact set $\{\tau_k\}_{k\ge 0}$ is $\bp^{\eta,a}_\go$--a.s. unbounded.
By the arguments given in \S\ref{sec:building}, to obtain the weak convergence
of the full polymer measure $\bP^{a}_{N,\go}$, as $N\to\infty$ along $[N]=\eta$,
toward a limit law $\bP^{\eta,a}_\go$, it remains to
determine the law of the sign $\gs_{\rho + 1}$
of the last (infinite) excursion (the notation has been introduced in
\S\ref{sec:building}).

We start with the free case. We want to show that $\bP_{N,\go}^\rf \big( S_N > 0 \big)$
has a limit as $N\to\infty$ along $[N]=\eta$.
By conditioning on the last zero before~$N$ we get
\begin{align*}
    \bP_{N,\go}^\rf \big( S_N > 0 \big) \;=\; \frac{1}{Z^\rf_{N,\go}} \,
    \sum_{k \ge 0} \sum_{\gamma \in \bbS} \sum_{n=0}^{N-1} M^{k*}_{[0],\gamma}(n) \,
    \bigg( \frac12 \sum_{t>N-n}K(t) \bigg) \,,
\end{align*}
where $K(\cdot)$ has been defined before \eqref{eq:Feller2} and
$M^{k*}$ denotes the convolution of the kernel $M$ with itself $k$ times,
the convolution between two kernels $F$ and $G$ being defined by
\begin{equation*}
    (F\ast G)_{\ga,\gb}(n) \;:=\; \sum_{m=1}^{n-1} \sum_{\gamma \in \bbS} F_{\ga,\gamma}(m)
    G_{\gamma,\gb}(n-m)\,.
\end{equation*}
Therefore, using \eqref{eq:Feller2} and \eqref{eq:as_Z_del} and recalling the
definition \eqref{eq:matrixB}, we obtain
\begin{equation} \label{eq:genauf}
    \exists\, \limtwo{N\to\infty}{[N]=\eta}
    \bP_{N,\go}^\rf \big( S_N > 0 \big) \;=\; \frac{c_K}{\Lambda^\rf_{[0],\eta}}
    \, \sum_{\gamma \in \bbS} (1-B^{-1})_{[0],\gamma} \,.
\end{equation}

Next we consider the constrained case, where we focus instead on
$\bP_{N,\go}^\rc \big( S_{\lfloor N/2 \rfloor} > 0 \big)$.
Conditioning on the last zero before and on the first zero after
$\lfloor N/2 \rfloor$,   we can write
\begin{align*}
    \bP_{N,\go}^\rc & \big( S_{\lfloor N/2 \rfloor} > 0 \big) \\
    & =\; \frac{1}{Z^\rc_{N,\go}} \,
    \sum_{k,h \ge 0} \sum_{\gamma,\zeta \in \bbS} \sum_{n=0}^{\lfloor N/2 \rfloor}
    \sum_{m=\lceil N/2 \rceil}^{N-1} M^{k*}_{[0],\gamma}(n) \,
    \bigg( \frac12 \, K(m-n) \, e^{\go^{(0)}_\zeta} \bigg)
    \, M^{h*}_{\zeta,[N]}(N-m) \,,
\end{align*}
and using again \eqref{eq:Feller2} and \eqref{eq:as_Z_del} we obtain
\begin{equation} \label{eq:genauc}
    \exists\, \limtwo{N\to\infty}{[N]=\eta}
    \bP_{N,\go}^\rc \big( S_{\lfloor N/2 \rfloor} > 0 \big) \;=\; \frac{c_K}{2 \Lambda^\rc_{[0],\eta}}
    \, \sum_{\gamma \in \bbS} (1-B^{-1})_{[0],\gamma} \,
    \sum_{\zeta \in \bbS} e^{\go^{(0)}_\zeta} (1-B^{-1})_{\zeta,\eta}\,.
\end{equation}

Now it is easy to check that \eqref{eq:genauf} and \eqref{eq:genauc} give exactly
the probability, under the infinite volume polymer measure $\bP^{\eta,a}_\go$,
that the sign $\gs_{\rho+1}$ of the last (infinite) excursion equals~$+1$, and this
completes the proof.\qed

\section{Non uniqueness of the infinite volume measure}
\label{sec:Gibbs}

We want to show that all infinite volume measures $\bP_{\go}^{\eta,a}$ appearing
in the strictly delocalized regime $\gd_\go < 1$, see Theorem~\ref{th:infvol} and
Section~\ref{sec:infvol}, are in reality superpositions of only two measures $\bQ^+_\go$
and $\bQ^-_\go$, that are {\sl extremal Gibbs measures} for our system.
We split the exposition in two parts:
\begin{itemize}
\item in \S\ref{sec:cont} we show, by purely
combinatorial arguments, that the law of the {\sl contact set} under $\bP_{\go}^{\eta,a}$
is a superposition of two basic laws $\bq^+_\go$ and $\bq^-_\go$;
\item in \S\ref{sec:extr} we show that $\bq^+_\go$ and $\bq^-_\go$ can be extended
to two laws $\bQ^+_\go$ and $\bQ^-_\go$ for the whole process $\{S_n\}_n$ which are
extremal Gibbs measure for our system.
\end{itemize}

\subsection{Decomposition of the contact set law}
\label{sec:cont}

Let $\bp^{\eta, a}_\go$ denote the law of the contact set $(\tau_k)_{k\ge 0}$ under the
infinite volume measure $\bP^{\eta, a}_\go$. As it has been shown in \S\ref{sec:tl_del},
under $\bp^{\eta, a}_\go$ the process $(\tau_k)_{k\ge 0}$ is a Markov renewal process
with semi--Markov kernel $\Gamma^{\eta, a}$, defined in \eqref{eq:GGG}. More explicitly,
for every $n\in\N$ and for all $0=:k_0<k_1< \cdots < k_n$ we have:
\begin{equation} \label{eq:alg}
\begin{split}
\bp^{\eta,a}_{\go} \big( \{k_1, \ldots, k_n \} \big) & \; =\;
\Gamma^{\eta, a}_{[0],[k_1]}(k_1) \, \cdots\, \Gamma^{\eta, a}_{[k_{n-1}],[k_n]}(k_n - k_{n-1})\\
& \;=\; M_{[0],[k_1]}(k_1) \, \cdots\, M_{[k_{n-1}],[k_n]}(k_n - k_{n-1})
\frac{\Lambda^{a}_{[k_n],\eta}}{\Lambda^{a}_{[0],\eta}}\,,
\end{split}
\end{equation}
where $\Lambda^{a}_{\ga,\gb}$ is defined in \eqref{eq:lambda} and the basic kernel
$M_{\ga,\gb}(n)$ has been introduced in \eqref{eq:matrix}.

\smallskip

To express the law $\bp^{\eta, a}_\go$ as a superposition we are going to exploit
the algebraic structure of \eqref{eq:alg}. However the steps are more transparent
if carried out in a general setting, and one is led to the following definition:
we introduce the set $\cC$ defined by
\begin{equation} \label{eq:defC}
    \cC \;:=\; \bigg\{ v \in (0,\infty)^\bbS \;:\
    \sum_{\gb \in \bbS} \bigg( \sum_{n\in\N} M_{\ga,\gb}(n) \bigg) v_\gb \le v_\ga\,,
    \quad \forall \ \ga\in\bbS \bigg\}\,.
\end{equation}
More explicitly, we recall that $B_{\ga,\gb}:= \sum_{n\in\N} M_{\ga,\gb}(n)$ has
spectral radius $\delta_\go<1$, and therefore we have
\begin{equation} \label{eq:defC2}
    \cC=\left\{v=(I-B)^{-1}w, \ w\in [0,\infty)^\bbS\backslash\{0\} \right\}\,.
\end{equation}
The reason for
such a definition is that if (and only if) $v\in\cC$ then the kernel $M_{\ga,\gb}(n) \cdot v_\gb/v_\ga$
is a (defective) semi--Markov kernel, that is
$\sum_{\gb, n}M_{\ga,\gb}(n)\, v_\gb/v_\ga \le 1$ for every $\ga \in\bbS$.
Therefore, for all $v\in\cC$, we can define a (defective) law $\bq^v$
for the contact set $\{\tau_k\}_{k\in\N}$ by
\begin{equation}
\label{eq:def_Q}
\begin{split}
    & \bq^v \big( \{k_1, \ldots, k_n \}\big)\; :=\; M_{[0],[k_1]}(k_1) \, \cdots\,
    M_{[k_{n-1}],[k_n]}(k_n - k_{n-1}) \cdot
    \frac{v_{[k_n]}}{v_{[0]}}\,,
\end{split}
\end{equation}
for every $n\in\N$ and for all $0=:k_0<k_1< \cdots < k_n$.

\smallskip

Now let us take two arbitrary vectors $v^+, v^- \in \cC$.
Since $\cC$ is a convex set, for all $p\in [0,1]$
the vector $v:= p v^+ + (1-p) v^-$ belongs to $\cC$,
hence the law $\bq^v$ is well-defined. The crucial result is expressed by
the following combinatorial lemma.

\smallskip

\begin{lemma} \label{lem:gibbs}
The law $\bq^{pv^+ + (1-p)v^-}$ is a superposition of the laws $\bq^{v^+}$ and $\bq^{v^-}$:
\begin{equation}\label{eq:Qq}
    \bq^{v} \;=\; r \,\bq^{v^+} + (1-r)\, \bq^{v^-}\,, \qquad
    \text{where} \qquad r = \frac{p v^+_{[0]}}{p v^+_{[0]} + (1-p) v^-_{[0]}} \in [0,1]\,.
\end{equation}
\end{lemma}

\medskip
\noindent
{\it Proof.}
By \eqref{eq:def_Q}, all we have to verify is that for every~$\ga \in \bbS$
\begin{equation}
    r \, \frac{v^+_{\ga}}{v^+_{[0]}} \ + \ (1-r) \,
    \frac{v^-_{\ga}}{v^-_{[0]}} \ = \
    \frac{p\,v^+_\ga\, + \, (1-p)\, v^-_\ga}
    {p\,v^+_{[0]}\, + (1-p)\, v^-_{[0]}}\,.
\end{equation}
By the definition \eqref{eq:Qq} of~$q$, we can rewrite the l.h.s. above as
\begin{equation}
\begin{split}
    & \frac{p\, v^+_{[0]}}{p\, v^+_{[0]} + (1-p)\, v^-_{[0]}}
    \, \frac{v^+_{\ga}}{v^+_{[0]}} \ + \
    \frac{(1-p)\, v^-_{[0]}}{p\, v^+_{[0]} + (1-p)\, v^-_{[0]}} \,
    \frac{v^-_{\ga}}{v^-_{[0]}} = \\
    & = \ \frac{p\, v^+_{\ga}}{p\, v^+_{[0]} + (1-p)\, v^-_{[0]}}
    \ + \ \frac{(1-p)\, v^-_{\ga}}{p\, v^+_{[0]} + (1-p)\, v^-_{[0]}}
    \ = \ \frac{p\,v^+_\ga\, + \, (1-p)\, v^-_\ga}
    {p\,v^+_{[0]}\, + (1-p)\, v^-_{[0]}}\,,
\end{split}
\end{equation}
and the proof is completed. \qed

\medskip

Next we come back to our model. We define two vectors $v^+(\go)$ and $v^-(\go)$ by
\begin{equation}\label{eq:vv}
v^+(\go)_\ga \, := \, \sum_{\gamma\in\bbS}
(1-B)^{-1}_{\ga,\gamma} \qquad  \qquad
v^-(\go)_\ga \, := \, \sum_{\gamma\in\bbS}
(1-B)^{-1}_{\ga,\gamma} \, e^{-\Sigma_{[0],\gamma}}
\end{equation}
where $B_{\ga,\gb} = B^\go_{\ga,\gb} := \sum_{n\in\N} M_{\ga,\gb}^\go(n)$,
see \eqref{eq:matrixB}, and $\Sigma_{\ga,\gb}$ is defined in \eqref{eq:def_Sigma}.
From \eqref{eq:defC2} we have that $v^\pm(\go) \in \cC$, and
the corresponding laws $\bp^{v^\pm(\go)}$ will be simply denoted by $\bq^\pm_\go$:
\begin{equation} \label{eq:defq+q-}
    \bq^+_\go := \bq^{v^+(\go)} \qquad \quad \bq^-_\go := \bq^{v^-(\go)}\,.
\end{equation}
We are ready to state the main result of this paragraph.

\medskip
\begin{proposition} \label{prop:decomp}
For every $a=\rf,\rc$ and $\eta \in \bbS$, the measures $\bp^{\eta,a}_\go$
are superpositions of the two laws $\bq^+_\go$ and $\bq^-_\go$:
\begin{equation} \label{eq:decomp}
    \bp^{\eta,a}_\go \;=\; r(\eta,a,\go) \, \bq^+_\go \;+\; \big(1-r(\eta,a,\go)\big) \, \bq^-_\go\,,
\end{equation}
with $r(\eta,a,\go) \in (0,1)$.
\end{proposition}

\medskip
\noindent{\it Proof.}
We introduce the vector
$v(\eta,a,\go)_\ga := \Lambda^a_{\ga, \eta}$ (the dependence of $\Lambda^a_{\ga,\gb}$
on $\go$ has not been explicitly indicated, but of course is present), and
notice that the law $\bp^{\eta,a}_\go$
coincides with $\bq^{v(\eta,a,\go)}$, cf. \eqref{eq:alg} and \eqref{eq:def_Q}.

To prove that $\bp^{\eta,a}_\go=\bq^{v(\eta,a,\go)}$ is a superposition of
$\bq^\pm_\go = \bq^{v^\pm(\go)}$, we are going
to exploit Lemma~\ref{lem:gibbs}. Let us be more precise: we are going to show that, for every
$a=\rf,\rc$ and $\eta \in \bbS$, the vector $v(\eta,a,\go)$ is a {\sl linear combination}
of two vectors $v^+(\go)$ and $v^-(\go)$ with positive coefficients:
\begin{equation}\label{eq:toprove}
    v(\eta,a,\go) \;=\; x \, v^+(\go) \;+\; y \, v^-(\go)\,, \qquad \quad
    x,y\in\R^+\,.
\end{equation}
Then the vector $v(\eta,a,\go)$ can be written as the following {\sl convex combination}:
\begin{equation} \label{eq:cv}
    v(\eta,a,\go) \;=\; \frac{x}{x+y} \, w^+(\go) \;+\; \frac{y}{x+y} \, w^-(\go)\,, \qquad \quad
    w^\pm(\go) := (x+y) \, v^\pm(\go)\,,
\end{equation}
and Lemma~\ref{lem:gibbs} yields that
$\bq^{v(\eta,a,\go)} = \bp^{\eta,a}_\go$ is a superposition of the two laws
$\bq^{w^\pm(\go)}$. However it is straightforward to see from \eqref{eq:def_Q}
that the laws $\bq^{w^\pm(\go)}$ are the same as $\bq^{v^\pm(\go)}$, because the vectors $v^\pm(\go)$
and $w^\pm(\go)$ differ only by a scale factor. Therefore from \eqref{eq:toprove}
it follows indeed that $\bp^{\eta,a}_\go$ is a superposition of $\bq^\pm_\go$, that is
what we have to prove.

Therefore it only remains to show that \eqref{eq:toprove}
holds true, where of course $x=x(\eta,a,\go)$
and $y=y(a,\eta,\go)$. We consider first the constrained case $a=\rc$: from the
definition \eqref{eq:lambda} of $\Lambda^\rc_{\ga,\eta}$ and from the definition
\eqref{eq:as_tPhi} of the matrix~$L$, we can write for $\ga \in\bbS$
\begin{align*}
    v(\rc,\eta,\go)_\ga &\; =\; \Lambda^\rc_{\ga,\eta}
    \; = \; \big[ (1-B)^{-1} L (1-B)^{-1} \big]_{\ga,\eta}\\
    &\; =\; \frac{c_K}{2} \sum_{\gamma,\zeta \in \bbS} (1-B)^{-1}_{\ga,\gamma} \big(1 +
    \exp(\Sigma_{\gamma, \zeta})\big) e^{\go^{(0)}_\zeta} (1-B)^{-1}_{\zeta,\eta}\,.
\end{align*}
Observing that $\Sigma_{\gamma,\zeta} = \Sigma_{[0],\zeta} - \Sigma_{[0],\gamma}$ and recalling
the definition \eqref{eq:vv} of $v^\pm(\go)$ we obtain
\begin{align*}
    v(\rc,\eta,\go)_\ga & = \Bigg( \frac{c_K}{2} \sum_{\zeta \in \bbS}
    e^{\go^{(0)}_\zeta} (1-B)^{-1}_{\zeta,\eta} \Bigg) v^+(\go)_\ga +
    \Bigg( \frac{c_K}{2} \sum_{\zeta \in \bbS}
    e^{\go^{(0)}_\zeta + \Sigma_{[0],\zeta}} (1-B)^{-1}_{\zeta,\eta} \Bigg) v^-(\go)_\ga\,,
\end{align*}
which shows that \eqref{eq:toprove} holds true for $a=\rc$ and gives an explicit
expression for $x(\rc,\eta,\go)$ and $y(\rc,\eta,\go)$.
With analogous (and simpler) arguments, for the free case we get
\begin{align*}
    v(\rf,\eta,\go)_\ga & \ = \ c_K \, v^+(\go)_\ga  \ + \
    \big( c_K e^{\Sigma_{[0],\eta}}\big) \, v^-(\go)_\ga\,.
\end{align*}
Thus \eqref{eq:toprove} holds true
also for $a=\rf$, with $x(\rf,\eta,\go) = c_K$ and $y(\rf,\eta,\go) = c_K e^{\Sigma_{[0],\eta}}$,
and the proof is completed.\qed

\medskip
Finally, we observe that one can obtain
an explicit formula for the weight $r(\eta,a,\go)$ appearing
in \eqref{eq:decomp}. From the expression for $r$ given in \eqref{eq:Qq}
and from \eqref{eq:cv} it follows that
\begin{align*}
    r(\eta,a,\go) & \;=\; \frac{x(\eta,a,\go)\, v^+(\go)_{[0]}}{x(\eta,a,\go)\, v^+(\go)_{[0]}
    \,+\, y(\eta,a,\go)\, v^-(\go)_{[0]}}\\
    & \;=\; \frac{x(\eta,a,\go)\, v^+(\go)_{[0]}}{v(\eta,a,\go)_{[0]}} \;=\;
    \frac{x(\eta,a,\go)\, v^+(\go)_{[0]}}{\Lambda^a_{[0],\eta}} \,,
\end{align*}
having used \eqref{eq:toprove} and the definition $v(\eta,a,\go)_\ga := \Lambda^a_{\ga,\eta}$.
Observe that
the precise values of $x(\eta,a,\go)$ is the coefficient of $v^+(\go)_\ga$ in the last two
equations of the proof of Proposition~\ref{prop:decomp}, cf. \eqref{eq:toprove}.
Then, recalling the definition \eqref{eq:vv} of $v^\pm(\go)$,
we obtain the following formula for $r(\eta,a,\go)$: for the constrained case $a=\rc$
\begin{align} \label{eq:supperc}
    r(\eta,\rc,\go) & \;=\; \frac{\sum_{\gamma\in\bbS} \, (1-B)^{-1}_{[0],\gamma} \cdot
    \frac{c_K}{2} \sum_{\zeta \in \bbS}
    e^{\go^{(0)}_\zeta} (1-B)^{-1}_{\zeta,\eta}}{\Lambda^\rc_{[0],\eta}}\,,
\end{align}
and for the free case $a=\rf$
\begin{equation}\label{eq:supperf}
    r(\eta, \rf,\go) \;=\; \frac{\sum_{\gamma\in\bbS} \, (1-B)^{-1}_{[0],\gamma}
    \cdot c_K}{\Lambda^\rf_{[0],\eta}}\,.
\end{equation}
The exact value of $r(\eta,a,\go)$ will be important in the next paragraph.

\subsection{Extremal Gibbs measures}
\label{sec:extr}

The aim of this paragraph is to show that the decomposition of the contact set law
$\bp^{\eta,a}_\go$ in terms of the two laws $\bq^\pm_\go$, proved in the previous paragraph,
can be lifted from the space of the contact set $\{\tau_n\}_n$
to the space of trajectories of $\{S_n\}_n$. More precisely, we are going to show that
for all $a=\rf,\rc$ and $\ga \in \bbS$ the infinite volume measure $\bP^{\eta,a}_\go$
is a superposition of two laws $\bQ^\pm_\go$, depending only on~$\go$, which have
$\bq^\pm_\go$ as contact set laws and which are extremal Gibbs measures for our system.

\smallskip

Let us first recall some basic notions. A measure $\bQ$ on $\Z^{\N\cup\{0\}}$ is said
to be a {\sl Gibbs measure} for our system if it satisfies the so--called DLR equation,
that in our setting reads as follows:
for all $M\in\N$ and for all $A\subset \Z^M$ we have
\begin{equation}\label{eq:dlr}
\bQ\big( (S_1,\ldots,S_M)\in A\ \big| \, S_M\, \big) \, = \,
\bP^\rf_{M,\go} \big( (S_1,\ldots,S_M) \in A\ \big| \, S_M\,\big) \qquad
\bQ\text{--a.s.}\,.
\end{equation}
The set of all Gibbs measures is clearly a convex set, that is
if $\bQ_1$ and $\bQ_2$ are Gibbs measure and $p\in[0,1]$ then the convex combination
$p\bQ_1 + (1-p)\bQ_2$ is a Gibbs measure too. If a Gibbs measure $\bQ$ cannot
be written as a nontrivial convex combination of two distinct Gibbs measures, then
$\bQ$ is said to be {\sl extremal}.
The standard reference on Gibbs measures is \cite{cf:Geo}.

Both the free and the constrained polymer measures $\bP^\rf_{N,\go}$ and $\bP^\rc_{N,\go}$
satisfy relation \eqref{eq:dlr} for any $M\le N$. Then it is not a surprise that any
weak limit of $\bP^a_{N,\go}$, as $N\to\infty$ along a subsequence,
satisfies \eqref{eq:dlr} for all $M\in\N$, that is it is a Gibbs measure, cf.~\cite[Th.~4.17]{cf:Geo}.
In particular, all infinite volume measures $\bP^{\eta,a}_\go$ for $a=\rf,\rc$ and
$\eta\in\bbS$, that are found in Theorem~\ref{th:infvol}, are Gibbs measures.

\smallskip

The basic Gibbs measures $\bQ^\pm_\go$ extending $\bq^\pm$ are introduced in the next lemma.

\smallskip
\begin{lemma} \label{lem:Q}
There exist two {\sl extremal} Gibbs measures $\bQ^+_\go$ and $\bQ^-_\go$
such that the law of the contact set $(\tau_n)_{n\ge 0}$ under $\bQ^\pm_\go$
is exactly $\bq^\pm_\go$. Moreover these laws satisfy
\begin{equation} \label{eq:asQ+Q-}
    \bQ^+_\go \Big( \lim_{N\to\infty} S_N = +\infty \Big) \;=\; 1 \qquad \quad
    \bQ^-_\go \Big( \lim_{N\to\infty} S_N = -\infty \Big) \;=\; 1\,.
\end{equation}
\end{lemma}
\medskip

\noindent
The proof of this lemma is given below. Now that we have introduced the two laws
$\bQ^\pm_\go$, we are ready to state and prove the main result of this section.

\smallskip
\begin{proposition}
For all $a=\rf,\rc$ and $\eta \in \bbS$, the infinite volume measures $\bP^{\eta,a}_\go$
given in Theorem~\ref{th:infvol}, for $\gd_\go < 1$, are superpositions of the two
laws $\bQ^+_\go$ and $\bQ^-_\go$ given in Lemma~\ref{lem:Q}. More precisely:
\begin{equation} \label{eq:dec}
    \bP^{\eta,a}_\go \;=\; r(\eta,a,\go) \, \bQ^+_\go \;+\;(1-r(\eta,a,\go)) \, \bQ^-_\go\,,
\end{equation}
where the weight $r(\eta,a,\go) \in (0,1)$ is given by \eqref{eq:supperc} and \eqref{eq:supperf}
for $a=\rc,\rf$ respectively.
\end{proposition}

\medskip\noindent
{\it Proof.}
We already know by Proposition~\ref{prop:decomp} that relation \eqref{eq:dec} holds
true if restricted to events involving only the contact set, see \eqref{eq:decomp}.
Now notice that,
conditionally on the level set, the law of the signs and of the moduli of the excursions
(except for the last infinite one) are {\sl the same} under the three laws
$\bP^{\eta,a}_\go$, $\bQ^+_\go$ and $\bQ^-_\go$, that is they are given by \eqref{defsigma}
and \eqref{defexc}: this is just because all three laws are Gibbs measures for our
system and hence satisfy the relation \eqref{eq:dlr}.
Therefore relation \eqref{eq:dec} holds true if restricted to the events that happen
not later than the last contact point (more precisely, restricted on the
$\gs$--field $\gs(\tau_\rho, S_k: 0 \le k \le \tau_\rho)$, where
$\rho := \sup\{k\ge 0: \tau_k < +\infty\}$ is the index of the last contact point).

Then it remains to focus on the sign $\gs_{\rho + 1}$ and on the modulus
$e_{\rho + 1}(\cdot)$ of the last (infinite) excursion (the notation has been introduced
in \S\ref{sec:prel}). For the modulus $e_{\rho + 1}(\cdot)$ there are no problems, because
it has the same law under each of $\bP^{\eta,a}_\go$, $\bQ^+_\go$ and $\bQ^-_\go$,
see \eqref{defexclul}. About the sign $\gs_{\rho + 1}$, we know from Lemma~\ref{lem:Q}
that under $\bQ^+_\go$ it is $+1$ and under $\bQ^-_\go$ it is $-1$, hence under the r.h.s. of
\eqref{eq:dec} the variable $\gs_{\rho + 1}$ takes the values $+1$ and $-1$ with probabilities
respectively equal to $r(\eta,a,\go)$ and $1-r(\eta,a,\go)$. However, the l.h.s. of \eqref{eq:dec},
that is $\bP^{\eta,a}_\go$, gives exactly the same law to $\gs_{\rho+1}$, cf. \eqref{eq:genauf}
and \eqref{eq:genauc} with \eqref{eq:supperf} and \eqref{eq:supperc}, and this completes the
proof.\qed

\bigskip\noindent
{\it Proof of Lemma~\ref{lem:Q}.}
Let us introduce two modified finite volume polymer measures $\bP^+_{N,\go}$ and
$\bP^-_{N,\go}$, defined by
\begin{equation} \label{eq:defP+P-}
    \frac{\dd\bP^+_{N,\go}}{\dd \bP} (S) \;:=\; \frac{\exp\big(\cH'_N(S)\big)}{Z_{N,\go}^+} \,
    \ind_{(S_N > 0)} \qquad \
    \frac{\dd\bP^-_{N,\go}}{\dd \bP} (S) \;:=\; \frac{\exp\big(\cH'_N(S)\big)}{Z_{N,\go}^-} \,
    \ind_{(S_N < 0)}\,,
\end{equation}
and notice that $Z^\pm_{N,\go} = Z_{N,\go}^\rf \cdot \bP^\rf_{N,\go}(S_N \gtrless 0)$, cf. \eqref{eq:newP}.
Then from Theorem~\ref{th:as_Z} and equation \eqref{eq:genauf} it follows
that for any fixed $k\ge 0$, as $N\to\infty$ along $[N]=\eta$
\begin{equation*}
    Z_{N-k,\theta_k\go}^+ \sim \Bigg(\sum_{\gamma\in\bbS}(1-B)^{-1}_{[k],\gamma}
    \Bigg) \frac{c_K}{\sqrt N} \qquad
    Z_{N-k,\theta_k\go}^- \sim \Bigg(\sum_{\gamma\in\bbS}(1-B)^{-1}_{[k],\gamma}
    e^{-\Sigma_{[0],\gamma}} \Bigg) \frac{c_K\, e^{\Sigma_{[0],\eta}}}{\sqrt N}\,.
\end{equation*}
Therefore for every fixed $k\ge 0$ we obtain
\begin{equation*}
    \exists \lim_{N\to\infty} \ \frac{Z^+_{N-k,\theta_k\go}}{Z^+_{N,\go}}
    \, = \, \frac{v^+(\go)_{[k]}}{v^+(\go)_{[0]}} \qquad \ \ \
    \exists \lim_{N\to\infty} \ \frac{Z^-_{N-k,\theta_k\go}}{Z^-_{N,\go}}
    \, = \, \frac{v^-(\go)_{[k]}}{v^-(\go)_{[0]}}\,,
\end{equation*}
where the vectors $v^\pm(\go)$ have been defined in \eqref{eq:vv}. But then, following
closely the proof of Proposition~\ref{pr:infvol2}, it is easy to prove that both
the measures $\bP^\pm_{N,\go}$ converge weakly on $\Z^{\N\cup\{0\}}$ as $N\to\infty$
toward two limit measures, that we denote by $\bQ^\pm_\go$, such that the contact set
$\{\tau_n\}_{n\ge 0}$ under $\bQ^\pm_\go$ has law $\bq^\pm_\go$, cf. \eqref{eq:defq+q-}.
In particular, the cardinality of the contact set $\{\tau_n\}_{n\ge 0}$ is $\bQ^\pm_\go$--a.s.
finite.
Moreover, by the definition \eqref{eq:defP+P-} of $\bP^\pm_{N,\go}$, it follows that
the sign of the last (infinite) excursion under $\bQ^+_\go$ (resp. under $\bQ^-_\go$)
is deterministic and takes the value $+1$ (resp. $-1$). This proves \eqref{eq:asQ+Q-}.

Being weak limit of finite volume polymer measures with suitable boundary conditions,
the two laws $\bQ^\pm$ are automatically Gibbs measures for our system, cf.~\cite[Th.~4.17]{cf:Geo}.
To complete the proof,
it only remains to show that they are extremal, and by \cite[Th.~7.7]{cf:Geo} this is
equivalent to showing that they are trivial on the tail $\gs$--field of the sequence
$\{S_n\}_{n\ge 0}$.

Let us denote by $\cG_n:= \gs ( S_k: k\ge n )$ the $\gs$--field generated by the
variables $\{S_k\}$ with index $k\ge n$.
We recall that the tail $\gs$--field $\cT$ is defined by
$\cT := \bigcap_{m \in \N} \cG_m$. Let us denote by $\Theta^{-1}$ the {\sl inverse
shift} defined on $\cG_1$, that is for $A\in\cG_1$ the event $\Theta A \in \cG_0$ is defined by
\begin{equation*}
    (S_0, S_1, S_2, \ldots) \in \Theta^{-1} A \qquad \iff \qquad (S_1, S_2, S_3, \ldots ) \in A\,.
\end{equation*}
By iteration we can define the $n$--shift $\Theta^{-n}$ on $\cG_n$, in particular
if $A \in \cT$ then $\Theta^{-n} A$ is well defined for all $n\in\N$ and
$\Theta^{-n} A\in \cT$.

We have to show that $\bQ^\pm_\go(A) = 0$ for all $A\in\cT$, and for conciseness we
focus on $\bQ^+_\go$ (the case $\bQ^-_\go$ is analogous). We recall that
$\rho := \sup\{k\ge 0: \tau_k < +\infty\}$ denotes
the index of the last contact point, and we stress that $\bQ^+_\go(\rho < +\infty)=1$.
We also recall from \S\ref{sec:building} that the last excursion $\{e_{\rho+1}(k)\}_{k \ge 0}
:= \{S_{\tau_\rho + k}\}_{k\ge 0}$ has under $\bQ^+_\go$ the law $\bP^\uparrow$ of the
{\sl random walk conditioned to stay positive}, see \eqref{defexclul} and \cite{cf:BerDon}.
We point out that $\bP^\uparrow$ is the law of a Markov chain on $\N\cup\{0\}$ which is
transient: $\bP^\uparrow (\lim_{N\to\infty} S_N = +\infty) = 1$, cf.~\cite{cf:BerDon}.

By conditioning on the value of the last contact point, we can write
\begin{equation} \label{eq:end}
    \bQ^+_\go (A ) \;=\; \sum_{n\ge 0} \bQ^+_\go( A \,|\, \tau_{\rho} = n ) \;
    \bQ^+_\go \big( \tau_{\rho} = n \big)\,.
\end{equation}
However if $A\in\cT$ then $A \in \cG_n$, for all~$n$, hence
\begin{equation*}
    \bQ^+_\go( A \,|\, \tau_{\rho} = n ) \;=\;
    \bQ^+_\go\big( \{S_{n+k}\}_{k\ge 0} \in A \,\big|\, \tau_{\rho} = n \big) \;=\;
    \bP^\uparrow \big( \Theta^{-n} A \big)\,.
\end{equation*}
We have already remarked that $\Theta^{-n} A\in \cT$ for all~$n$, hence if we
show that the law $\bP^\uparrow$ is trivial on~$\cT$ then from \eqref{eq:end}
it follows that $\bQ^+_\go(A)=0$ and this completes the proof.

\smallskip

Let $k\in\N$, $\ell_1,\ldots,\ell_k\in\N$ and set $M_n:=\bP^\uparrow(S_i=\ell_i,
i=1,\ldots,k \, | \, \cG_n)$, $n>k$. Then $(M_n)_{n>k}$ is a $(\cG_n)_{n>k}$--inverse
martingale, hence $M_n$ converges $\bP^\uparrow$--a.s. and in $L^1(\dd\bP^\uparrow)$
to $M:=\bP^\uparrow(S_i=\ell_i, i=1,\ldots,k \, | \, \cT)$. On the other hand, by the Markov property:
\begin{align} \label{eq:bau}
\begin{split}
\bP^\uparrow(S_i=\ell_i, i=1,\ldots,k \, | \, \cG_n) &\; = \;
\bP^\uparrow(S_i=\ell_i, i=1,\ldots,k \, | \, S_n) \\
&\;=\;\left[ \prod_{i=1}^k p^\uparrow_1(\ell_{i-1},\ell_i)
\right] \frac{p^\uparrow_{n-k}(\ell_k,S_n)}{p^\uparrow_n(0,S_n)}\,,
\qquad \bP^\uparrow\text{--a.s.}\,,
\end{split}
\end{align}
where $l_0:=0$ and $p_j^\uparrow(a,b)$ is the $j$--th iteration of the
transition kernel of $\bP^\uparrow$, that
is the $\bP^\uparrow$--probability that
$S$ goes from $a$ to $b$ in $j$ steps. Now we claim that for every $x\in\N$
\begin{equation}\label{eq:?}
\lim_{n\to\infty} \frac{p^\uparrow_{n-k}(x,S_n)}{p^\uparrow_n(0,S_n)}
\; = \; 1\,, \qquad \bP^\uparrow\text{--a.s.}\,.
\end{equation}
Then we obtain:
\[
\bP^\uparrow(S_i=\ell_i, i=1,\ldots,k \, | \, \cT) \, = \,
\left[ \prod_{i=1}^k p^\uparrow_1(\ell_{i-1},\ell_i)
\right] \, = \, \bP^\uparrow(S_i=\ell_i, i=1,\ldots,k)
\]
and it follows that $\cT$ is independent of $\gs(S_i:\,i=1,\ldots,k)$.
Since this is true for all $k\in\N$, $\cT$ is independent of itself
and therefore must be trivial.

It remains to prove \eqref{eq:?}. We recall that $\bP^\uparrow$ is the law
of a transient Markov chain on~$\R$, cf.~\cite{cf:BerDon}, and we denote by
$\bP^\uparrow_x$ the law with starting point~$x\in\N\cup\{0\}$. Then we can rephrase
\eqref{eq:?} in the following way: for all $k, x \in \N$
\begin{equation}\label{eq:??}
    \lim_{n\to\infty} \frac{\bP^\uparrow_x(S_{n-k}=y)}{\bP^\uparrow_0(S_n = y)} \, \bigg|_{y=S_n}
    \;=\; 1\,,\qquad \bP^\uparrow\text{--a.s.}\,.
\end{equation}
We stress that we already know that the l.h.s. of this equation
has a limit as $n\to\infty$, $\bP^\uparrow$--a.s. (it suffices to give a look
at the r.h.s. of \eqref{eq:bau} and to recall that the l.h.s. of \eqref{eq:bau}
converge $\bP^\uparrow$--a.s. by the martingale argument outlined above).
Therefore it suffices to show that, for $\bP^\uparrow$--a.e. $S=\{S_n\}_n$, there exists
a subsequence $(n_k)_k = ( n_k(S) )_k$ such that the l.h.s. of \eqref{eq:??} tends to~$1$ as
$n\to\infty$ along the subsequence $(n_k)_k$.

Let us denote by $(\{R_t\}_{t\ge 0},P)$ a standard Bessel(3) process starting from zero.
Then from \cite[Th.~5.1]{cf:BryDon} we have that, for any fixed~$x$, under $\bP^\uparrow_x$
the sequence $S_n / (\sqrt{2p}\sqrt n)$ convergence in law toward $R_1$ (note that $\sqrt{2p}$
is the variance of $S_1$ under the unperturbed random walk measure $\bP$). In particular,
for all $0<a<b<\infty$ we have
\begin{equation} \label{eq:rel0}
    \bP^\uparrow_x \Big( S_n \in \Big[a \sqrt{2p}\, \sqrt n \,,\, b \sqrt{2p}\, \sqrt n\Big] \Big)
    \;\to\; P\big( R_1 \in \big[a, b \big]\big) \qquad \quad (n\to\infty)\,.
\end{equation}
It follows that, for $\bP_x^\uparrow$--a.e. $S=\{S_n\}_n$, there exists
a subsequence $(n_k)_k = ( n_k(S) )_k$ such that $\liminf_k S_{n_k} / \sqrt{n_k} > 0$ and
$\limsup_k S_{n_k} / \sqrt{n_k} < \infty\,$: indeed, from \eqref{eq:rel0} and
by Fatou's lemma
\begin{equation*}
    \bP^\uparrow_x \Big( S_n \in \Big[a \sqrt{2p}\, \sqrt n \,,\, b \sqrt{2p}\, \sqrt n\Big]
    \ \text{i.o.} \, \Big) \;\ge\; P\big( R_1 \in \big[a, b \big]\big)\,,
\end{equation*}
where $\{A_n \text{ i.o.}\} := \limsup_n A_n$ for a sequence of events $(A_n)_n$, and
since $P( R_1 \in [a, b ])$ can be made arbitrarily close to $1$ by choosing $a$ small and $b$ large,
we obtain the claim.

Therefore, in order to prove \eqref{eq:??} it is enough to show that for any
sequence $(y_n)_n\subset\N$ such that $\liminf_n y_n/ \sqrt n > 0$ and $\limsup_n y_n/ \sqrt n<\infty$ we have
\begin{equation}\label{i.e.}
    \frac{\bP^\uparrow_x(S_{n-k}=y_n)}{\bP^\uparrow_0(S_n = y_n)} \;\to\; 1 \qquad \quad (n\to\infty)\,.
\end{equation}

We denote by $p(\cdot,\cdot)$ the transition kernel of $S$ under $\bP$:
\[
p(x,y) \, := \, p \, \ind_{(y=x+1)} + p \, \ind_{(y=x-1)} + (1-2p) \, \ind_{(y=x)}, \qquad x,y\in\N\,,
\]
and we denote by $\bP_x$ the law of $x+S$ under $\bP$, $x\in\N\cup\{0\}$.
We recall that the transition probability kernel $p^\uparrow(x,y)
=p_1^\uparrow(x,y)$ of $S$ under $\bP^\uparrow$
is a $h$-transform of $p(x,y) \ind_{(y>0)}$:
\[
p^\uparrow(x,y) \, = \, p(x,y) \, \ind_{(y>0)} \, \frac{h(y)}{h(x)}, \qquad x,y\geq 0\,,
\]
where $h:\N\cup\{0\}\mapsto(0,\infty)$ satisfies:
\[
\sum_{y} p(x,y) \, \ind_{(y>0)} \, h(y) \, = \, h(x), \qquad x\geq 0\,,
\]
see \cite{cf:BerDon}. It is easy to see that necessarily $h(x)=h(0) \, x/p$ for
all $x\geq 1$, hence for all $x,y\geq 1$:
\begin{align*}
    \frac{\bP^\uparrow_x(S_{n-k}=y)}{\bP^\uparrow_0(S_n = y)} & \; = \;
\frac{\bP_x(S_n=y, S_1>0, ..., S_n>0)}{\bP_0(S_n=y, S_1>0, ..., S_n>0)} \cdot
\frac{p}{x}\\
    & \;=\; \frac{\bP_x(S_n=y)- \bP_x(S_n=-y)}{p\,[\bP_0(S_{n-1}=y-1)-\bP_0(S_{n-1}=y+1)]} \cdot
\frac{p}{x}\,,
\end{align*}
where we have used the reflection principle.
The Local Limit Theorem given by \cite[Th.~13 in
Ch.~VII.3]{cf:petrov} yields the expansion:
\[
\bP_0(S_n=z\cdot{\sqrt n}) \, = \, \gga(z)
\left[1+\frac c{\sqrt n}\, (z^3-3z)\right] + o(1/{\sqrt n})\,,
\]
uniformly in $z\in n^{-1/2}\, \N$, where $\gga(\cdot)$ is the density of $\cN(0,1)$ and $c$
a positive constant. Then we obtain:
\[
\frac{\bP_x(S_n=y_n)- \bP_x(S_n=-y_n)}{\bP_0(S_{n-1}=y_n-1)-\bP_0(S_{n-1}=y_n+1)} \;\to\;
x\,, \qquad \quad (n\to\infty)\,,
\]
and the proof of \eqref{i.e.} is complete.\qed

\subsection*{Acknowledgments}
G.G. acknowledges the support of GIP-ANR, project {\sl POLINTBIO}.

\end{document}